\documentclass[11pt]{article}

\usepackage[a4paper,margin=1in]{geometry}
\usepackage[T1]{fontenc}
\usepackage[utf8]{inputenc}
\usepackage{lmodern}
\usepackage{microtype}
\usepackage{xcolor}
\usepackage{amsmath,amssymb,amsfonts,mathtools}
\usepackage{amsthm}
\usepackage{bm}
\usepackage{booktabs}
\usepackage{multirow}
\usepackage{hyperref}
\usepackage{graphicx}
\usepackage{subcaption}
\usepackage{algorithm}
\usepackage{algorithmic}
\usepackage{url}
\usepackage[backend=biber,style=numeric,sorting=none]{biblatex}
\usepackage{hyperref}
\usepackage[nameinlink,capitalize,noabbrev]{cleveref}
\usepackage{enumitem}
\usepackage{float}
\usepackage{placeins}
\usepackage{multicol}

\hypersetup{colorlinks=true,linkcolor=blue,citecolor=blue,urlcolor=blue}
\setlist[itemize]{topsep=0.3em,itemsep=0.2em}

\newtheorem{theorem}{Theorem}[section]
\newtheorem{proposition}[theorem]{Proposition}
\newtheorem{lemma}[theorem]{Lemma}
\newtheorem{corollary}[theorem]{Corollary}
\theoremstyle{remark}

\newcommand{\R}{\mathbb{R}}
\newcommand{\E}{\mathbb{E}}
\newcommand{\F}{\mathcal{F}}
\newcommand{\T}{\mathcal{T}}
\newcommand{\diag}{\operatorname{diag}}
\newcommand{\tr}{\operatorname{tr}}

\title{NAG-GS: A Semi-Implicit Stochastic Optimizer with Large-Step Stability}
\author{
Valentin Leplat\thanks{Corresponding author. Department of Mathematics and Operational
Research, Universit\'e de Mons, Rue de Houdain 9, 7000 Mons, Belgium
(\href{mailto:Valentin.LEPLAT@umons.ac.be}{\texttt{Valentin.LEPLAT@umons.ac.be}}).
VL acknowledges support from the European Research Council
(ERC Consolidator Grant, eLinoR, No.~101085607).}
\and Daniil Merkulov
\and Aleksandr Katrutsa
\and Daniel Bershatsky
\and Olga Tsymboi
\and Ivan Oseledets
}
\date{}

\makeatletter
\newcommand{\compactcontents}{%
  \begingroup
  \begin{minipage}{\linewidth}%
  \setcounter{tocdepth}{1}%
  \hypersetup{linkcolor=black}%
  \footnotesize\setlength{\parskip}{0pt}%
  \setlength{\multicolsep}{4pt}%
  \renewcommand{\l@section}{\@dottedtocline{1}{0em}{1.5em}}%
  \textbf{Contents}\par
  \begin{multicols}{2}\@starttoc{toc}\end{multicols}%
  \end{minipage}%
  \endgroup
}
\makeatother
\begin{document}
\maketitle

\begin{abstract}
Stochastic first-order methods are widely used for large-scale optimization and machine learning, where their performance often depends on the learning rate.
We study NAG-GS, a stochastic momentum method based on a semi-implicit Gauss--Seidel discretization of a Nesterov-type flow.
For strongly convex quadratics, we determine the exact range of stable step sizes for fixed damping and the step size that minimizes the worst-case spectral radius.
When the damping equals the strong-convexity parameter, the method recovers the accelerated asymptotic factor $(\sqrt{\kappa}-1)/(\sqrt{\kappa}+1)$, where $\kappa$ is the condition number.
We also prove that the quadratic iteration converges linearly with dynamic damping when its limiting fixed-damping iteration is stable.
For general strongly convex, possibly nonsmooth objectives, we establish an $O(1/k)$ bound on the expected objective gap of the last iterate, using a conditionally unbiased stochastic first-order oracle with bounded second moment.
The result holds with dynamic damping from any $\gamma_0>0$.
We show that NAG-GS can also be written as a Heavy-Ball iteration with explicitly determined step sizes and momentum coefficients.
We test the stochastic result on nonsmooth hinge-loss classification problems, and further evaluate NAG-GS on quadratic problems, logistic regression, convolutional neural networks, RoBERTa, and Vision Transformers.
The experiments show that NAG-GS can remain effective at learning rates where standard optimizers become unstable or lose accuracy.
\end{abstract}


 \section{Introduction}
\label{sec:intro}

Stochastic gradient methods are the standard tools for training large machine-learning models. 
In practice, momentum is often added to improve the behavior of stochastic gradient descent (SGD) during the first iterations, and adaptive methods such as AdamW are widely used when the learning rate is difficult to choose. 
This makes stability with respect to the learning rate an important practical issue.

A useful way to design first-order methods is to start from a continuous-time model and then choose a discretization. 
Different discretizations of the same differential equation can have very different stability properties. 
This point has been studied for gradient flows and accelerated flows; see, for example,~\cite{su2014differential,shi2019acceleration,zhang2018direct,luo2021differential}. 
In particular, Luo and Chen~\cite{luo2021differential} considered a Nesterov-type flow and a semi-implicit Gauss--Seidel discretization. 
The deterministic state equations used in NAG-GS coincide with the Gauss--Seidel discretization of Luo and Chen. {Starting from this discretization, we study its exact stability on strongly convex quadratics, prove stochastic convergence with dynamic damping, and evaluate NAG-GS on large-scale learning problems.}

In this paper, we study the resulting method, called NAG-GS, in deterministic and stochastic settings. 
We consider two choices for the damping parameter. 
The dynamic rule follows directly from the ODE discretization. 
The fixed rule with \(\gamma=\mu\) is an exact special case of the dynamic rule. 
More generally, keeping \(\gamma\) fixed gives a useful variant that can be tuned independently. 
It also makes it possible to derive an exact stability analysis on quadratic problems. 
We keep both choices throughout the paper and state explicitly which one is used in each result.

We focus on three questions. 
First, how large can the step size be on strongly convex quadratics, and which stable step gives the fastest asymptotic contraction? 
Second, what can be proved when the damping parameter is updated dynamically? 
Third, can one obtain a stochastic convergence result for general strongly convex functions under the usual conditionally unbiased first-order oracle? 
These questions lead to complementary results: an exact constant-step analysis on quadratics, a convergence result for the dynamic damping sequence, and an $O(1/k)$ stochastic result with decreasing steps.

\subsection{Related work}

\paragraph{ODE discretization and NAG-GS.}
Continuous-time interpretations of accelerated methods have been developed in many works; see, among others, \cite{su2014differential,muehlebach2019dynamical,shi2021understanding,zhang2018direct}. Luo and Chen \cite{luo2021differential} derive the Nesterov accelerated gradient flow used here and study several discretizations. Their Gauss--Seidel scheme has exactly the deterministic state equations that form the core of NAG-GS. They also explain that a spectral bound by itself does not give a general Lyapunov contraction for nonlinear objectives. Our quadratic analysis stays in the linear setting, where the spectral calculation can be made exact. We also study the stochastic recursion and the dynamic damping rule.

\paragraph{Lyapunov analysis and quasi-monotone methods.}
Wilson, Recht and Jordan \cite{JMLR:v22:20-195} develop a Lyapunov view of accelerated and momentum methods. Their quasi-monotone method becomes algebraically identical to NAG-GS in Euclidean geometry when the damping is fixed to \(\gamma=\mu\). This connection gives an \(O(1/k)\) deterministic rate for strongly convex, possibly nonsmooth functions under bounded subgradients.   For dynamic damping, we use the time-varying energy of Luo and Chen~\cite[Lemma~3]{luo2021differential}. We extend their estimate to a conditionally unbiased stochastic oracle.

\paragraph{Stochastic acceleration.}
Lan's AC-SA method \cite{lan2012optimal} gives optimal complexity guarantees for smooth stochastic convex and composite problems. Its theory is more general than the theory developed here. Laborde and Oberman \cite{laborde2020lyapunov} also study stochastic accelerated dynamics through Lyapunov arguments. Our focus is different: we study the stability of the semi-implicit NAG-GS update, derive exact quadratic step-size results, and test the method in large nonconvex learning problems. We include AC-SA in the convex numerical comparisons.

\paragraph{Practical momentum methods.}
SGD with momentum and AdamW are standard choices in deep learning \cite{loshchilov2019decoupled}. NAG-GS has one auxiliary state vector, as does SGD with momentum, while Adam-like methods keep first- and second-moment vectors. Our numerical study compares these methods on convex and nonconvex problems. The nonconvex tests are empirical: the theory in this paper does not claim to describe the full behavior of deep neural networks.

\subsection{Contributions}

Our main contributions are the following.
\begin{enumerate}[label=\arabic*.]
    \item For strongly convex, possibly nonsmooth objectives, we prove an $O(1/k)$ expected last-iterate rate with dynamic damping from any $\gamma_0>0$. We assume that
    \[
       \E[g_{k+1}\mid\F_k]\in\partial f(x_{k+1}),
       \qquad
       \E[\|g_{k+1}\|^2\mid\F_k]\leq G^2.
    \]
    \item  We give a direct modal analysis of NAG-GS on strongly convex quadratics. For fixed damping, we obtain the exact Schur-stability interval. For \(\gamma\geq\mu\), we also obtain the step size that minimizes the worst-case spectral radius. This separates two quantities that play different roles: the largest stable step and the step giving the fastest asymptotic contraction.
    \item For \(\gamma=\mu\), the optimal asymptotic factor is \((\sqrt{\kappa}-1)/(\sqrt{\kappa}+1)\). Thus NAG-GS recovers the classical accelerated factor on strongly convex quadratics.
    \item We analyze the dynamic damping rule. We first give its exact expression for general positive step sizes. With a constant step size, \(\gamma_k\) converges geometrically to \(\mu\). We then prove that the time-varying quadratic iteration is linearly convergent whenever its limiting fixed-damping matrix is Schur stable.
    \item We show that NAG-GS can be written as a Heavy-Ball iteration. Its step sizes and momentum coefficients are determined by the NAG-GS parameters, and the two forms produce the same iterates when they use the same initialization and oracle directions. This holds for fixed and dynamic damping.
    \item {We test the theoretical results on quadratic and nonsmooth strongly convex problems, check the agreement with the matched Heavy-Ball iteration, and evaluate NAG-GS on larger machine-learning problems.}
 \end{enumerate}

The rest of the paper is organized as follows. \Cref{sec:method} recalls the continuous-time model and defines the NAG-GS family. \Cref{sec:quadratic} gives the quadratic stability analysis, including the dynamic damping result. \Cref{sec:stochastic} gives the stochastic \(O(1/k)\) result. \Cref{sec:experiments} presents the numerical experiments. Additional proofs and numerical details are collected in the appendices.

\section{Derivation and formulation of NAG-GS}
\label{sec:method}

\subsection{Continuous-time model}

Accelerated first-order methods can be interpreted through continuous-time systems involving the optimization variable together with an auxiliary momentum variable. This viewpoint is useful here because different discretizations of the same continuous dynamics can lead to very different numerical stability properties.

We use the Nesterov-type system considered by Luo and Chen~\cite{luo2021differential} as our starting point. Here \(x(t)\) is the optimization variable, \(v(t)\) is an auxiliary momentum variable, and \(\gamma(t)>0\) controls the damping of the dynamics. For a differentiable objective \(f:\mathbb{R}^n\to\mathbb{R}\), the system is

\begin{equation}
\label{eq:flow}
\begin{aligned}
    \dot x(t) &= v(t)-x(t),\\
    \dot v(t) &= \frac{\mu}{\gamma(t)}\bigl(x(t)-v(t)\bigr)-\frac{1}{\gamma(t)}\nabla f(x(t)),\\
    \dot\gamma(t) &= \mu-\gamma(t),\qquad \gamma(0)=\gamma_0>0.
\end{aligned}
\end{equation}
For a quadratic objective \(f(x)=\tfrac12 x^\top A x\), the first two equations can be written as a linear system. 
The parameter \(\gamma(t)\) satisfies
\[
    \gamma(t)=\mu+(\gamma_0-\mu)e^{-t}.
\]
Thus \(\gamma(t)\) approaches \(\mu\), and it is constant when \(\gamma_0=\mu\).

In a stochastic setting, the exact gradient is replaced by noisy first-order information.
At the diffusion level, suppose that the gradient contribution is perturbed by a Brownian term.
Since the gradient appears in the second equation multiplied by \(1/\gamma_t\), a gradient-noise
coefficient \(\Sigma_t\) would give a diffusion term of the form
\(-\gamma_t^{-1}\Sigma_t\,\mathrm dW_t\) in the dynamics of \(v_t\).
For simplicity, we denote the resulting effective diffusion coefficient by \(\sigma\) and consider
the additive-noise model
\begin{equation}
\label{eq:sde}
\begin{aligned}
    \mathrm d x_t &= (v_t-x_t)\,\mathrm dt,\\
    \mathrm d v_t
    &= \left[
        \frac{\mu}{\gamma_t}(x_t-v_t)
        -\frac{1}{\gamma_t}\nabla f(x_t)
    \right]\mathrm dt
    + \sigma\,\mathrm dW_t.
\end{aligned}
\end{equation}
Here \(W_t\) is a standard Brownian motion, and \(\sigma\) represents the effective noise level
in the \(v_t\)-dynamics. We use \eqref{eq:sde} only as a simple continuous-time motivation for
the stochastic discretization. Similar diffusion viewpoints are common in stochastic optimization;
see, for example,~\cite{malladi2022sdes}.

The stochastic convergence result in \Cref{sec:stochastic} is entirely discrete and does not rely
on this diffusion model. It uses a standard conditionally unbiased first-order oracle with bounded
second moment and, in particular, does not assume Gaussian gradient noise.

\subsection{Why a semi-implicit discretization?}
\label{subsec:why-gs}

The choice of the discretization is important here. A fully explicit Euler step is cheap, but its stable step-size range can be small. A fully implicit step has better stability properties, but for a general objective it would require solving an implicit equation involving the gradient at every iteration. This is not attractive for large-scale problems.

The Gauss--Seidel update lies between these two choices. To see this, consider first the quadratic case \(f(x)=\tfrac12x^\top A x\) and keep \(\gamma>0\) fixed. With \(y=(x,v)\), the first two equations of \eqref{eq:flow} have the form \(\dot y=G_\gamma y\), where
\[
G_\gamma=
\begin{bmatrix}
-I & I\\
(\mu I-A)/\gamma & -(\mu/\gamma)I
\end{bmatrix}.
\]
Split this matrix as \(G_\gamma=M_\gamma+N_\gamma\), with
\[
M_\gamma=
\begin{bmatrix}
-I & 0\\
(\mu I-A)/\gamma & -(\mu/\gamma)I
\end{bmatrix},
\qquad
N_\gamma=
\begin{bmatrix}
0 & I\\
0 & 0
\end{bmatrix}.
\]
A Gauss--Seidel time step is then
\begin{equation}
\label{eq:gs-splitting}
    \frac{y_{k+1}-y_k}{\alpha_k}
    =M_{\gamma_k}y_{k+1}+N_{\gamma_k}y_k.
\end{equation}
The new value \(x_{k+1}\) is therefore used immediately in the update of \(v_{k+1}\). At the same time, the scheme does not require the solution of a linear system involving \(A\). For a general function, the same construction is obtained by replacing \(Ax_{k+1}\) by a first-order direction evaluated at \(x_{k+1}\).

This semi-implicit structure changes the stability properties of the discretization. As we show in \Cref{sec:quadratic}, it allows relatively large stable steps on some quadratic problems. It is also computationally simple: \(x_{k+1}\) is computed first, one gradient or stochastic-gradient call is made at this point, and \(v_{k+1}\) is then obtained by vector operations only.

\subsection{Semi-implicit Gauss--Seidel update}

With step size \(\alpha_k>0\), the state equations are
\begin{equation}
\label{eq:gs-discretization}
\begin{aligned}
    \frac{x_{k+1}-x_k}{\alpha_k} &= v_k-x_{k+1},\\
    \frac{v_{k+1}-v_k}{\alpha_k} &= \frac{\mu}{\gamma_k}(x_{k+1}-v_{k+1})-\frac{1}{\gamma_k}g_{k+1},
\end{aligned}
\end{equation}
where \(g_{k+1}\) is a deterministic or stochastic first-order direction evaluated at \(x_{k+1}\). Solving the two implicit relations gives the more convenient form
\begin{equation}
\label{eq:naggs-state}
\boxed{
\begin{aligned}
    x_{k+1} &= \frac{x_k+\alpha_k v_k}{1+\alpha_k},\\
    v_{k+1} &= \frac{\gamma_k v_k+\alpha_k\mu x_{k+1}-\alpha_k g_{k+1}}{\gamma_k+\alpha_k\mu}.
\end{aligned}}
\end{equation}
This expression remains meaningful for \(\mu=0\) as long as \(\gamma_k>0\).

There are two damping choices that we will use.

\paragraph{Dynamic damping.}
Applying an implicit step to \(\dot\gamma=\mu-\gamma\) gives
\[
    \frac{\gamma_{k+1}-\gamma_k}{\alpha_k}=\mu-\gamma_{k+1}.
\]
Solving this scalar relation gives
\begin{equation}
\label{eq:gamma-dynamic}
    \gamma_{k+1}=\frac{\gamma_k+\alpha_k\mu}{1+\alpha_k}.
\end{equation}
When \(\gamma_0=\mu\), this gives \(\gamma_k\equiv\mu\).

\paragraph{Fixed damping.}
For analysis and for practical tuning, we also consider \(\gamma_k\equiv\gamma>0\). The fixed choice \(\gamma=\mu\) is an exact special case of the dynamic rule. The fixed regime is not meant to replace the dynamic one. It isolates the state dynamics from the scalar damping recursion and makes it possible to determine the stability boundary exactly.

In the convex theory, \(\mu\) denotes the strong-convexity parameter of the objective. In the nonconvex experiments, there is no such known quantity. We therefore treat the parameters denoted by \(\mu\) and \(\gamma\) simply as algorithmic parameters controlling the state update and tune them empirically. The update remains exactly \eqref{eq:naggs-state}, but these parameters are not interpreted as the strong-convexity and damping parameters of the continuous-time theory. We only use parameter choices for which the denominator in \eqref{eq:naggs-state} is nonzero; unstable trials are discarded.

\begin{algorithm}[t]
\caption{NAG-GS}
\label{alg:naggs}
\begin{algorithmic}[1]
\REQUIRE Initial point \(x_0\), set \(v_0=x_0\), step sizes \(\alpha_k>0\), parameter \(\mu\), and positive damping \(\gamma_0\).
\FOR{$k=0,1,2,\ldots$}
    \STATE $x_{k+1}\leftarrow (x_k+\alpha_k v_k)/(1+\alpha_k)$.
    \STATE Query a first-order direction $g_{k+1}$ at $x_{k+1}$.
    \STATE $v_{k+1}\leftarrow (\gamma_k v_k+\alpha_k\mu x_{k+1}-\alpha_k g_{k+1})/(\gamma_k+\alpha_k\mu)$.
    \STATE Either update $\gamma_{k+1}\leftarrow(\gamma_k+\alpha_k\mu)/(1+\alpha_k)$, or keep a fixed damping $\gamma_{k+1}\leftarrow\gamma_k$.
\ENDFOR
\end{algorithmic}
\end{algorithm}

The dominant cost of NAG-GS is the oracle call. The method stores one auxiliary vector \(v_k\), so its state size is comparable to SGD with momentum and smaller than Adam-like methods that keep two moment vectors.

\subsection{Relation to existing methods}

The deterministic state equations \eqref{eq:gs-discretization} coincide with the Gauss--Seidel discretization studied by Luo and Chen~\cite{luo2021differential}. Here, we consider the corresponding stochastic update and study its stability and convergence properties.

A second useful connection appears when \(\gamma_k\equiv\mu\). Then \eqref{eq:naggs-state} becomes
\begin{equation}
\label{eq:qmm-euclidean}
\begin{aligned}
    x_{k+1} &= \frac{x_k+\alpha_k v_k}{1+\alpha_k},\\
    v_{k+1} &= \frac{v_k+\alpha_k x_{k+1}-(\alpha_k/\mu)g_{k+1}}{1+\alpha_k}.
\end{aligned}
\end{equation}
This is the Euclidean quasi-monotone update of Wilson, Recht and Jordan \cite{JMLR:v22:20-195}, after identifying their step parameter with \(\alpha_k\). We use this relation in \Cref{sec:stochastic}.

 \paragraph{Heavy-Ball representation.}
The auxiliary variable $v_k$ can be removed from the NAG-GS update. Indeed, the first equation of \eqref{eq:naggs-state} gives
\[
    v_k=x_{k+1}+\frac{x_{k+1}-x_k}{\alpha_k}.
\]
Substituting this expression into the second equation of \eqref{eq:naggs-state} gives
\begin{equation}
\label{eq:hb-minimal}
    x_{k+2}=x_{k+1}-\eta_k g_{k+1}+\beta_k(x_{k+1}-x_k),
\end{equation}
where
\begin{equation}
\label{eq:hb-coeff-minimal}
\begin{aligned}
    \eta_k&=\frac{\alpha_{k+1}\alpha_k}{(1+\alpha_{k+1})(\gamma_k+\alpha_k\mu)},\\
    \beta_k&=\frac{\alpha_{k+1}\gamma_k}{\alpha_k(1+\alpha_{k+1})(\gamma_k+\alpha_k\mu)}.
\end{aligned}
\end{equation}
This is a Heavy-Ball iteration~\cite{polyak1964some} with step sizes and momentum coefficients determined by the NAG-GS parameters. If $v_0=x_0$, then $x_1=x_0$. Starting from these same two points and using the same oracle directions, the two formulations generate exactly the same iterates. This is an algebraic identity and holds for both fixed and dynamic damping.

In particular, $\alpha_k$ is not the gradient step size $\eta_k$ in \eqref{eq:hb-minimal}. Comparisons with Heavy Ball using a fixed momentum coefficient therefore compare different parameter choices, not different families of iterations.

\subsection{Scope of the analysis}
\label{subsec:scope}

The same state update is used in several parameter regimes, but the available results are different. \Cref{tab:theory-scope} summarizes what is proved in this paper. This distinction is important because the fixed and dynamic damping rules answer different questions. The fixed regime gives exact formulas on quadratics. The dynamic regime follows the ODE discretization and is treated as a time-varying linear system. The general stochastic result uses decreasing steps and   dynamic damping from any $\gamma_0>0$, including $\gamma_k\equiv\mu$ .

\begin{table}[ht!]
\centering
\small
\caption{Parameter regimes considered in the analysis.}
\label{tab:theory-scope}
\begin{tabular}{p{0.25\textwidth}p{0.16\textwidth}p{0.13\textwidth}p{0.30\textwidth}}
\toprule
Problem & Damping & Step size & Result \\
\midrule
Strongly convex quadratic, deterministic
& fixed \(\gamma>0\)
& constant
& exact Schur-stability interval; optimal fixed step for \(\gamma\geq\mu\) \\
\addlinespace
Strongly convex quadratic, deterministic
& dynamic \(\gamma_k\)
& constant
& linear convergence when the limiting \(\gamma=\mu\) system is Schur stable \\
\addlinespace
Strongly convex quadratic, stochastic
& fixed \(\gamma>0\)
& constant
& stationary covariance inside the Schur-stable region \\
\addlinespace
General strongly convex, possibly nonsmooth
&   dynamic $\gamma_k$
 & decreasing
& \(O(1/k)\) rate in expectation under a conditionally unbiased oracle \\
\addlinespace
Nonconvex learning problems
& fixed, tuned
& fixed/tuned experimentally
& empirical study only \\
\bottomrule
\end{tabular}
\end{table}

\section{Quadratic stability and acceleration}
\label{sec:quadratic}

We now study
\begin{equation}
\label{eq:quadratic}
    f(x)=\frac{1}{2} x^\top A x,
    \qquad A=A^\top\succ0,
\end{equation}
with
\[
    0<\mu=\lambda_{\min}(A)\leq \lambda_{\max}(A)=L<\infty.
\]
A linear term can be removed by shifting the minimizer, so \eqref{eq:quadratic} is sufficient for the analysis.

\subsection{Modal iteration}

Let \(u\) be an eigenvector of \(A\) with eigenvalue \(\lambda\in[\mu,L]\). For constant \(\alpha>0\) and fixed \(\gamma>0\), define
\[
    \xi_k=\langle x_k,u\rangle,
    \qquad
    \nu_k=\langle v_k,u\rangle.
\]
Then
\begin{equation}
\label{eq:modal}
\begin{bmatrix}\xi_{k+1}\\ \nu_{k+1}\end{bmatrix}
=T_\lambda(\alpha,\gamma)
\begin{bmatrix}\xi_k\\ \nu_k\end{bmatrix},
\end{equation}
where
\begin{equation}
\label{eq:Tlambda}
T_\lambda(\alpha,\gamma)=
\begin{bmatrix}
\dfrac{1}{1+\alpha} & \dfrac{\alpha}{1+\alpha}\\[0.3cm]
\dfrac{\alpha(\mu-\lambda)}{(1+\alpha)(\gamma+\alpha\mu)}
&
\dfrac{\gamma}{\gamma+\alpha\mu}
+\dfrac{\alpha^2(\mu-\lambda)}{(1+\alpha)(\gamma+\alpha\mu)}
\end{bmatrix}.
\end{equation}
Thus the full \(2n\)-dimensional iteration is an orthogonal direct sum of \(2\times2\) systems. This gives the \(n\)-dimensional result directly, without pairing the eigenvalues two by two.

The determinant and trace are
\begin{equation}
\label{eq:trace-det}
\det T_\lambda
=\frac{\gamma}{(1+\alpha)(\gamma+\alpha\mu)},
\end{equation}
\begin{equation}
\label{eq:trace}
\tr T_\lambda
=\frac{2\gamma+\alpha(\gamma+\mu)+\alpha^2(\mu-\lambda)}{(1+\alpha)(\gamma+\alpha\mu)}.
\end{equation}

\subsection{Exact Schur-stability interval}

A matrix is Schur stable if all its eigenvalues have modulus strictly smaller than one. In the quadratic case, this is exactly the condition needed for the deterministic linear iteration to converge to zero.

\begin{theorem}[Exact fixed-damping stability]
\label{thm:schur}
Let \(f\) be given by \eqref{eq:quadratic}, and let \(\alpha>0\) and \(\gamma>0\) be constant. Then NAG-GS is Schur stable for every eigendirection of \(A\) if and only if
\[
4\gamma+2\alpha(\gamma+\mu)+\alpha^2(2\mu-L)>0.
\]
Equivalently:
\begin{itemize}
    \item if \(L\leq2\mu\), every \(\alpha>0\) is stable;
    \item if \(L>2\mu\), stability holds exactly for \(0<\alpha<\alpha_{\rm stab}(\gamma)\), where
    \begin{equation}
    \label{eq:alpha-stab}
        \alpha_{\rm stab}(\gamma)
        =\frac{\gamma+\mu+\sqrt{\gamma^2-6\gamma\mu+\mu^2+4\gamma L}}{L-2\mu}.
    \end{equation}
\end{itemize}
At \(\alpha=\alpha_{\rm stab}(\gamma)\), the mode \(\lambda=L\) has an eigenvalue equal to \(-1\).
\end{theorem}

The proof is short. It applies the Jury conditions to the characteristic polynomial of \(T_\lambda\); see \Cref{app:fixed-proof}. The condition is worst in the eigendirection \(\lambda=L\), which is why only \(L\) appears in the final stability boundary.

\paragraph{How to read the stability result.}
The two cases in \Cref{thm:schur} have a simple interpretation. If \(L\leq2\mu\), or equivalently \(\kappa\leq2\), there is no finite upper stability bound on \(\alpha\). This does not mean that arbitrarily large steps give the fastest convergence. It only means that the iterates do not lose Schur stability as the step increases. When \(L>2\mu\), there is a finite boundary, and the mode associated with the largest eigenvalue is the first one to become unstable.

For the important choice \(\gamma=\mu\), these expressions can be written only in terms of the condition number \(\kappa=L/\mu\). If \(\kappa>2\),
\begin{equation}
\label{eq:alpha-stab-kappa}
    \alpha_{\rm stab}(\mu)
    =\frac{2}{\sqrt{\kappa-1}-1}.
\end{equation}
The best step obtained below is
\[
    \alpha_*=\frac{2}{\sqrt\kappa-1}.
\]
Hence \(\alpha_*<\alpha_{\rm stab}(\mu)\) whenever \(\kappa>2\). For large \(\kappa\), the two values are close. In other words, on ill-conditioned quadratics the fastest fixed step lies near the edge of the stability region. For moderate condition numbers, the gap can be more visible.

The formulas also become simpler after scaling the damping by \(\mu\). Let
\[
    \delta=\frac{\gamma}{\mu},
    \qquad
    \kappa=\frac{L}{\mu}.
\]
For \(\kappa>2\), the stability boundary is
\begin{equation}
\label{eq:alpha-stab-dimensionless}
    \alpha_{\rm stab}(\delta,\kappa)
    =\frac{\delta+1+\sqrt{\delta^2-6\delta+1+4\delta\kappa}}{\kappa-2}.
\end{equation}
For \(\delta\geq1\) and \(\kappa>1\), the optimal fixed step from \Cref{thm:optimal-fixed} is
\begin{equation}
\label{eq:alpha-opt-dimensionless}
    \alpha_*(\delta,\kappa)
    =\frac{\delta+1+\sqrt{(\delta-1)^2+4\delta\kappa}}{\kappa-1}.
\end{equation}
Thus, after normalization, the quadratic behavior depends only on the condition number and on the relative damping \(\gamma/\mu\). This is useful when comparing problems with different scales.

\subsection{Best fixed step size}

We next assume \(\gamma\geq\mu\). This includes the important case \(\gamma=\mu\).
Define
\begin{equation}
\label{eq:alpha-opt}
    \alpha_*(\gamma)
    =\frac{\gamma+\mu+\sqrt{(\gamma-\mu)^2+4\gamma L}}{L-\mu},
    \qquad L>\mu,
\end{equation}
and
\begin{equation}
\label{eq:r-gamma}
    r_\gamma(\alpha)=\frac{\gamma}{\gamma+\alpha\mu}.
\end{equation}

\begin{theorem}[Optimal fixed step size]
\label{thm:optimal-fixed}
Assume \(L>\mu>0\) and \(\gamma\geq\mu\). Let
\[
    R(\alpha,\gamma)=\max_{\lambda\in[\mu,L]}\rho\bigl(T_\lambda(\alpha,\gamma)\bigr).
\]
Then \(\alpha_*(\gamma)\) in \eqref{eq:alpha-opt} is the unique minimizer of \(R(\alpha,\gamma)\) over \(\alpha>0\). Moreover,
\begin{equation}
\label{eq:rho-opt-general}
    R\bigl(\alpha_*(\gamma),\gamma\bigr)
    =\frac{\gamma}{\gamma+\alpha_*(\gamma)\mu}.
\end{equation}
For every \(0<\alpha\leq\alpha_*(\gamma)\), all modal eigenvalues satisfy
\[
    |z|\leq r_\gamma(\alpha),
\]
and equality is attained in the mode \(\lambda=\mu\).
\end{theorem}

A proof is given in \Cref{app:fixed-proof}. The main point is again a two-dimensional Jury test, now applied after scaling the eigenvalues by \(r_\gamma(\alpha)\). For \(0<\alpha\leq\alpha_*(\gamma)\), the slowest mode is the one associated with \(\lambda=\mu\), and its spectral radius decreases as the step grows. At \(\alpha=\alpha_*(\gamma)\), the mode associated with \(\lambda=L\) reaches the same modulus. Increasing the step further makes the worst-case spectral radius grow again. Thus the largest stable step and the step giving the fastest asymptotic contraction are different quantities.

When \(\gamma=\mu\), \eqref{eq:alpha-opt} becomes
\begin{equation}
\label{eq:alpha-acc}
    \alpha_* =\frac{2\mu+2\sqrt{\mu L}}{L-\mu}
    =\frac{2}{\sqrt{\kappa}-1},
    \qquad \kappa=\frac{L}{\mu},
\end{equation}
and \eqref{eq:rho-opt-general} gives the classical accelerated factor.

\begin{corollary}[Accelerated quadratic factor]
\label{cor:accelerated-factor}
Under the assumptions of \Cref{thm:optimal-fixed}, set \(\gamma=\mu\). Then the optimal worst-case spectral radius is
\begin{equation}
\label{eq:accelerated-factor}
    R(\alpha_*,\mu)
    =\frac{\sqrt{\kappa}-1}{\sqrt{\kappa}+1}.
\end{equation}
For comparison, optimally tuned gradient descent has factor
\[
    \frac{\kappa-1}{\kappa+1}.
\]
\end{corollary}

The word ``factor'' in \Cref{cor:accelerated-factor} refers to the asymptotic spectral factor. In the mode \(\lambda=\mu\), the iteration matrix can be a Jordan block, so the Euclidean norm can contain a polynomial prefactor even though the asymptotic rate is governed by \eqref{eq:accelerated-factor}.

It is also useful to distinguish the case \(L=\mu\). Then every eigendirection is identical and, for \(\gamma=\mu\), the spectral radius is \(1/(1+\alpha)\). There is no finite optimal step: the asymptotic factor decreases to zero as \(\alpha\to\infty\). The formula \eqref{eq:alpha-acc} is therefore stated only for \(L>\mu\).

\subsection{Dynamic damping}
\label{subsec:dynamic}

We now return to the damping update \eqref{eq:gamma-dynamic}. We first record a simple property that is useful independently of the quadratic analysis.

\begin{proposition}[Evolution of the damping parameter]
\label{prop:gamma-sequence}
Let \(\gamma_0>0\) and \(\alpha_k>0\). Then
\begin{equation}
\label{eq:gamma-product}
    \gamma_k-\mu
    =(\gamma_0-\mu)
    \prod_{j=0}^{k-1}\frac{1}{1+\alpha_j}.
\end{equation}
In particular, \(\gamma_k\) stays between \(\gamma_0\) and \(\mu\). If \(\sum_{k=0}^{\infty}\alpha_k=+\infty\), then \(\gamma_k\to\mu\).
\end{proposition}

The identity follows by subtracting \(\mu\) from \eqref{eq:gamma-dynamic} and iterating the resulting scalar recursion. The convergence statement follows because the product in \eqref{eq:gamma-product} tends to zero whenever \(\sum_k\alpha_k=+\infty\). Thus the damping moves monotonically toward \(\mu\) for the usual positive step-size sequences.

For the quadratic convergence result below, we keep a constant step \(\alpha>0\). Then \eqref{eq:gamma-product} becomes
\begin{equation}
\label{eq:gamma-exact}
    \gamma_k-\mu
    =\frac{\gamma_0-\mu}{(1+\alpha)^k}.
\end{equation}
Hence the iteration matrix varies with \(k\), but it approaches the fixed-damping matrix with \(\gamma=\mu\) geometrically. We therefore obtain a time-varying linear system. In general, checking that every matrix in a time-varying product has spectral radius below one is not enough to conclude convergence. Here we use a different argument. The matrices converge to a Schur-stable limit. Therefore, after a finite number of iterations, they are all sufficiently close to the limiting matrix to be contractive in the same matrix norm. The first finitely many iterations only change the constant in the final convergence bound.

\begin{theorem}[Dynamic damping on quadratics]
\label{thm:dynamic}
Let \(f\) be given by \eqref{eq:quadratic}, let \(\gamma_0>0\), and use a constant step size \(\alpha>0\) together with the dynamic rule \eqref{eq:gamma-dynamic}. Assume that the limiting fixed-damping iteration with \(\gamma=\mu\) is Schur stable. Equivalently,
\begin{itemize}
    \item if \(L\leq2\mu\), take any \(\alpha>0\);
    \item if \(L>2\mu\), assume
    \begin{equation}
    \label{eq:dynamic-stab}
        0<\alpha<\frac{2\mu+2\sqrt{\mu(L-\mu)}}{L-2\mu}.
    \end{equation}
\end{itemize}
Then the deterministic dynamic-damping NAG-GS iterates converge linearly to the minimizer. More precisely, if \(\T_*\) denotes the full fixed-damping iteration matrix at \(\gamma=\mu\), then for every
\[
    q\in\bigl(\rho(\T_*),1\bigr)
\]
there exists \(C_q>0\) such that, for all sufficiently large \(k\),
\[
    \left\|\begin{bmatrix}x_k\\v_k\end{bmatrix}\right\|
    \leq C_q q^k
    \left\|\begin{bmatrix}x_0\\v_0\end{bmatrix}\right\|.
\]
\end{theorem}

The theorem does not require the matrices at different iterations to commute, and it does not require \(\gamma_0=\mu\). A finite number of early iterations can therefore lie far from the limiting damping value. What matters is that \(\gamma_k\) approaches \(\mu\) and that the limiting matrix is strictly inside the Schur-stable region. The proof uses \eqref{eq:gamma-exact}, continuity of the iteration matrix in \(\gamma\), and an induced norm in which the limiting matrix is contractive; see \Cref{app:dynamic-proof}.

This result is deliberately limited to deterministic quadratics with a constant step. It does not give the exact finite-iteration contraction factor of the time-varying system, and it does not give a general nonlinear convergence theorem for the dynamic damping rule. It does show, however, that the dynamic update is compatible with the stable fixed regime to which it converges.

In particular, choosing the optimal step of the limiting fixed-damping system, \(\alpha=\alpha_*\) from \eqref{eq:alpha-acc}, gives a dynamic scheme whose limiting spectral factor is the accelerated value \eqref{eq:accelerated-factor}.

\begin{corollary}[Asymptotic factor of the dynamic scheme]
\label{cor:dynamic-factor}
Under the assumptions of \Cref{thm:dynamic}, let \(z_k=(x_k,v_k)\). Then
\[
    \limsup_{k\to\infty}\|z_k\|^{1/k}
    \leq \rho(\T_*).
\]
If \(L>\mu\) and \(\alpha=\alpha_*\) from \eqref{eq:alpha-acc}, then
\[
    \limsup_{k\to\infty}\|z_k\|^{1/k}
    \leq
    \frac{\sqrt\kappa-1}{\sqrt\kappa+1}.
\]
\end{corollary}

This follows directly from \Cref{thm:dynamic}: the bound holds for every \(q>\rho(\T_*)\), and we then let \(q\) decrease to \(\rho(\T_*)\). The corollary concerns the asymptotic factor. It does not say that every finite iteration contracts by this amount.

\subsection{Additive noise and stationary covariance}
\label{subsec:stationary-cov}

For a fixed \(\alpha\) and \(\gamma\), the quadratic stochastic discretization can be written as
\begin{equation}
\label{eq:linear-stochastic}
    z_{k+1}=\T z_k+\xi_k,
\end{equation}
where \(z_k=(x_k,v_k)\) and \(\T\) is the full iteration matrix, orthogonally similar to the direct sum of the blocks in \eqref{eq:Tlambda}. If the Brownian model \eqref{eq:sde} is discretized with increments \(\Delta W_k\sim\mathcal N(0,\alpha I)\), the random part is
\[
    \xi_k=
    \begin{bmatrix}
    0\\[1mm]
    \dfrac{\sigma\gamma}{\gamma+\alpha\mu}\,\Delta W_k
    \end{bmatrix},
\]
so that
\begin{equation}
\label{eq:Q-brownian}
    Q:=\E[\xi_k\xi_k^\top]
    =
    \begin{bmatrix}
    0 & 0\\
    0 & \dfrac{\alpha\sigma^2\gamma^2}{(\gamma+\alpha\mu)^2}I
    \end{bmatrix}.
\end{equation}
The noise level does not change the Schur boundary of the deterministic matrix \(\T\). It does, however, determine how widely the stationary iterates are spread once the deterministic recursion is stable.

When \(\rho(\T)<1\), the mean converges to zero and the covariance converges to the unique solution of the discrete Lyapunov equation
\begin{equation}
\label{eq:lyapunov-cov}
    \Sigma_\infty=\T\Sigma_\infty\T^\top+Q.
\end{equation}
Equivalently,
\begin{equation}
\label{eq:vec-cov}
    \operatorname{vec}(\Sigma_\infty)
    =\bigl(I-\T\otimes\T\bigr)^{-1}\operatorname{vec}(Q).
\end{equation}
If the noise is Gaussian, the limiting stationary law is Gaussian. This gives two different quantities to consider when choosing the step: the transient contraction governed by \(\rho(\T)\) and the stationary spread governed by \eqref{eq:lyapunov-cov}. A step that is stable need not minimize either quantity. Near a finite Schur boundary, \(I-\T\otimes\T\) becomes nearly singular. The solution of the Lyapunov equation can therefore become very large, which corresponds to a large stationary spread of the iterates.

We only use \eqref{eq:lyapunov-cov} inside the stable region. Outside this region, the formal inverse in \eqref{eq:vec-cov} should not be interpreted as a stationary covariance. Additional illustrations are given in \Cref{app:quadratic-extra}.

\section{  Stochastic \texorpdfstring{$O(1/k)$}{O(1/k)}  convergence}
\label{sec:stochastic}

The quadratic analysis above concerns constant parameters. We now consider general strongly convex objectives and decreasing steps.   We allow dynamic damping from any $\gamma_0>0$; $\gamma_0=\mu$ recovers the  Euclidean quasi-monotone method \eqref{eq:qmm-euclidean}.

\subsection{Stochastic oracle}
\label{subsec:oracle}

Let \(f:\R^n\to\R\) be proper and \(\mu\)-strongly convex, with \(\mu>0\). We allow \(f\) to be nonsmooth and assume that it admits a minimizer \(x^*\); strong convexity makes this minimizer unique. We take the initialization \(x_0=v_0\) to be deterministic. Let \((\F_k)_{k\geq0}\) be the filtration generated by the algorithm up to iteration \(k\). The order of the NAG-GS update is important: \(x_{k+1}\) is computed from \(x_k\) and \(v_k\) before the oracle is queried. Hence \(x_{k+1}\) is \(\F_k\)-measurable, while \(g_{k+1}\) contains the new randomness of iteration \(k+1\).

We use the following stochastic oracle assumption:
\begin{equation}
\label{eq:oracle}
    s_{k+1}:=\E[g_{k+1}\mid\F_k]\in\partial f(x_{k+1}),
    \qquad
    \E[\|g_{k+1}\|^2\mid\F_k]\leq G^2.
\end{equation}
At each iteration, the oracle returns a random first-order direction at \(x_{k+1}\). Its conditional mean is a subgradient of \(f\) at this point, and its conditional second moment is bounded by \(G^2\) along the iterates.

In the proof, strong convexity is applied to
\[
    s_{k+1}=\E[g_{k+1}\mid\F_k],
\]
which is a valid subgradient of \(f\) at \(x_{k+1}\), while the random vector \(g_{k+1}\) is used in the actual update of \(v_{k+1}\). The difference \(g_{k+1}-s_{k+1}\) has zero conditional mean. Hence the corresponding cross term vanishes when we condition on \(\F_k\).

The second-moment bound in \eqref{eq:oracle} is only required along the iterates. For example, if the conditional variance satisfies
\[
    \E[\|g_{k+1}-s_{k+1}\|^2\mid\F_k]\leq\sigma_g^2
\]
and \(\|s_{k+1}\|\leq B\) along the iterates, then \eqref{eq:oracle} holds with \(G^2=B^2+\sigma_g^2\).

\subsection{Lyapunov estimate and rate}
\label{subsec:stoch-rate}

Define
\begin{equation}
\label{eq:energy}
    \mathcal E_k
    =A_k\left(f(x_k)-f(x^*)+  \frac{\gamma_k}{2} \|v_k-x^*\|^2\right),
\end{equation}
where \(x^*\) is the minimizer and the weights start from \(A_0>0\) and satisfy
\begin{equation}
\label{eq:A-recurrence}
    A_{k+1}=(1+\alpha_k)A_k.
\end{equation}

\begin{lemma}[One-step stochastic Lyapunov bound]
\label{lem:stoch-one-step}
Assume \eqref{eq:oracle}, use   the dynamic rule \eqref{eq:gamma-dynamic} , and let \(A_k\) satisfy \eqref{eq:A-recurrence}. Then
\begin{equation}
\label{eq:conditional-energy}
    \E[\mathcal E_{k+1}\mid\F_k]-\mathcal E_k
    \leq
      \frac{A_k\alpha_k^2}{2\gamma_k}
     \E[\|g_{k+1}\|^2\mid\F_k].
\end{equation}
\end{lemma}

The proof is given in \Cref{app:stochastic-proof}.   For deterministic gradients, this is Luo and Chen's estimate \cite[Lemma~3, (82)]{luo2021differential}
.

\begin{theorem}[Stochastic \(O(1/k)\) rate]
\label{thm:stoch-rate}
Assume \eqref{eq:oracle}. Use NAG-GS with   dynamic damping from any $\gamma_0>0$, and set $\underline\gamma=\min\{\gamma_0,\mu\}$,
 \begin{equation}
\label{eq:stoch-schedule}
    \alpha_k=\frac{2}{k+2},
    \qquad
    A_k=\frac{(k+2)(k+3)}{6}.
\end{equation}
Then, for every   \(k\geq1\),
 \begin{equation}
\label{eq:stoch-rate-explicit}
  \E[f(x_{k})-f(x^*)]
 \leq
  \frac{6\mathcal E_0}{(k+2)(k+3)}
 +
  \frac{3kG^2}{\underline\gamma(k+2)(k+3)} .
\end{equation}
In particular,
\[
      \E[f(x_{k})-f(x^*)] =O(1/  k ).
\]
\end{theorem}

The schedule \eqref{eq:stoch-schedule} satisfies \eqref{eq:A-recurrence}. More precisely,
\[
    A_k\alpha_k^2
    =\frac{2(k+3)}{3(k+2)}\leq1.
\]
 The damping recursion gives $\gamma_k\ge\underline\gamma>0$.  Taking total expectation in \eqref{eq:conditional-energy} and summing from   \(j=0\) to \(k-1\) gives
 \[
      \E[\mathcal E_{k}]
     \leq \mathcal E_0+  \frac{kG^2}{2\underline\gamma} .
\]
Since   \(\mathcal E_{k}\geq A_{k}(f(x_{k})-f(x^*))\) , division by   \(A_{k}=(k+2)(k+3)/6\)  gives \eqref{eq:stoch-rate-explicit}.   The fixed case $\gamma_k\equiv\mu$ is recovered when $\gamma_0=\mu$.

The bound in \Cref{thm:stoch-rate} concerns the objective value at $x_k$ itself, without averaging the iterates. It applies to general strongly convex, possibly nonsmooth objectives. For a deterministic subgradient oracle, the same argument gives the corresponding deterministic result. Unlike the constant-step quadratic analysis in \Cref{sec:quadratic}, this result uses decreasing steps and gives an $O(1/k)$ rate in expectation.

 \section{Experiments}
\label{sec:experiments}

The numerical study covers strongly convex quadratics, a nonsmooth strongly convex classification problem, multiclass logistic regression, VGG-11, ResNet-20, RoBERTa and Vision Transformers. The code for the large-scale experiments is available at \url{https://github.com/naggsopt/naggs}.

The experiments have four roles. The quadratic tests can be compared directly with the spectral analysis and show how the admissible step-size range changes with the condition number. The {nonsmooth tests illustrate convergence and the effect of dynamic damping.} Logistic regression then gives a smooth convex but nonquadratic problem with the fixed learning rates used in practice. Finally, the neural-network and Transformer experiments test whether the ability to use large learning rates remains useful in nonconvex training.

The theory in \Cref{sec:quadratic,sec:stochastic} does not cover the nonconvex deep-learning experiments. In those tests, \(\gamma\) is kept fixed and tuned as a damping hyperparameter. The parameter denoted by \(\mu\) in the implementation is also treated as an algorithmic parameter and is not interpreted as an unknown strong-convexity modulus. We keep the notation \(\mu\) below to match the implementation. The purpose of these experiments is therefore empirical; they are not presented as a verification of the strongly convex theory.

\subsection{Strongly convex quadratic problem}
\label{sec:exp-quadratic}

We first consider
\[
    f(x)=\frac{1}{2}x^\top A x-b^\top x,
\]
where \(A\in\R^{100\times100}\) is symmetric positive definite, \(\mu=\lambda_{\min}(A)\), and \(L=\lambda_{\max}(A)\). We declare convergence when \(f(x_k)-f^*\leq10^{-4}\). Divergent runs are shown at the sentinel value \(10^{10}\) in the figure; this value is used only to mark divergence. NAG-GS uses the fixed choice \(\gamma=\mu\). We compare it with gradient descent (GD), accelerated gradient descent (AGD), Heavy Ball (HB), and AC-SA \cite{lan2012optimal}. For HB we use momentum \(0.9\). For AC-SA, we vary the step-size-scale parameter used in our implementation. We test 70 learning rates on a logarithmic grid over \([10^{-3},10]\).

\begin{figure}[ht!]
    \centering
    \begin{subfigure}{0.48\linewidth}
        \includegraphics[width=\linewidth]{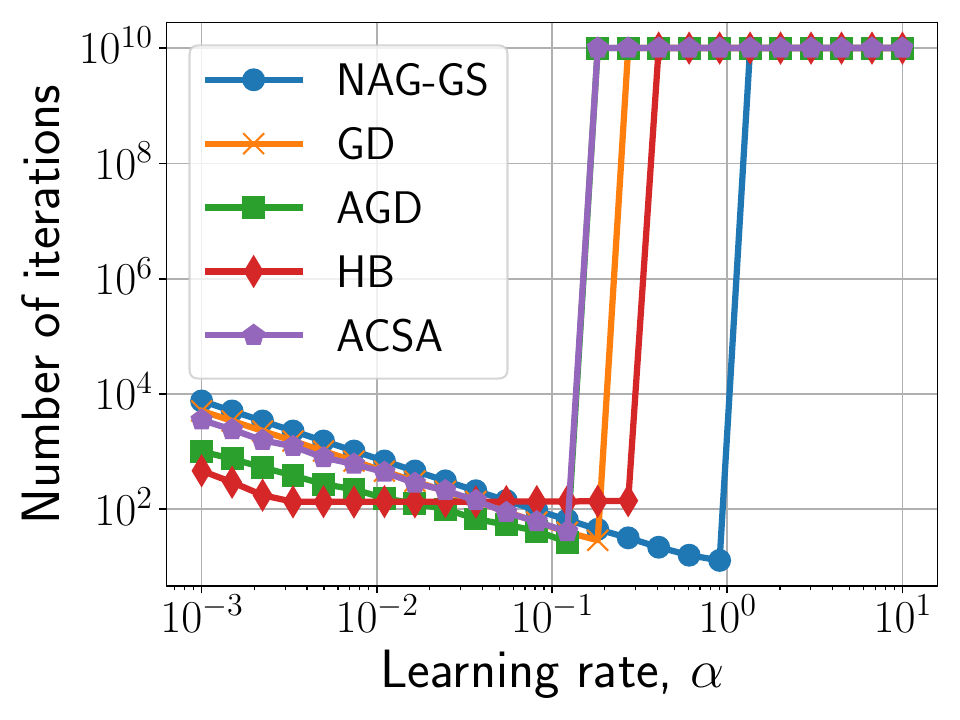}
        \caption{$\mu=1$, $L=10$.}
    \end{subfigure}
    \hfill
    \begin{subfigure}{0.48\linewidth}
        \includegraphics[width=\linewidth]{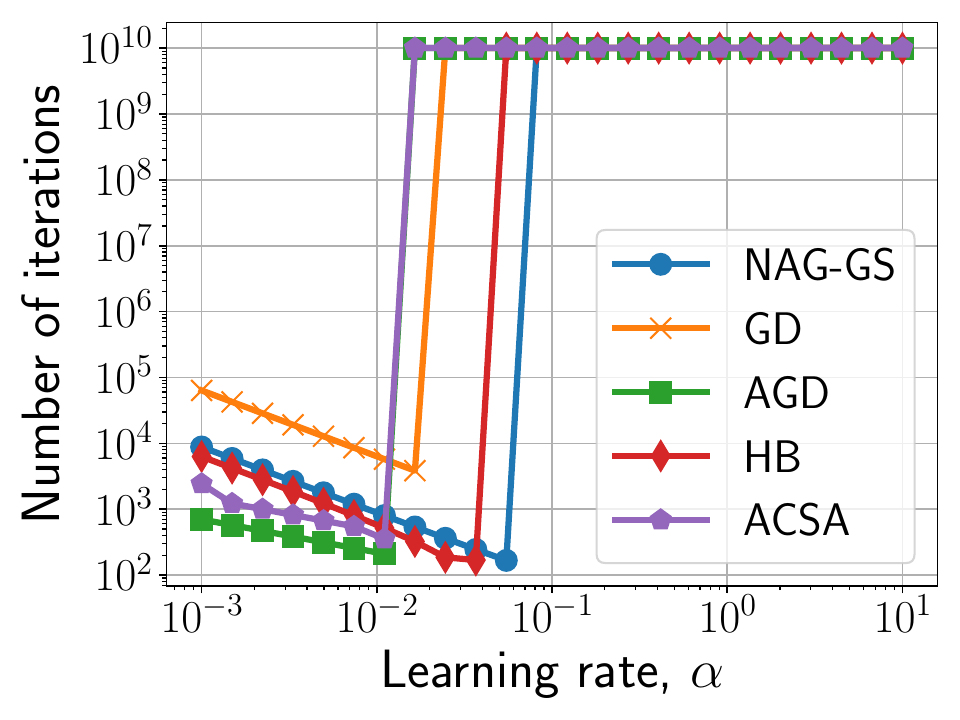}
        \caption{$\mu=10^{-1}$, $L=100$.}
    \end{subfigure}
    \caption{Number of iterations needed to reach the target tolerance as a function of the learning rate. The value $10^{10}$ indicates divergence within the test protocol.}
    \label{fig:quadratic-lr}
\end{figure}

\Cref{fig:quadratic-lr} shows two condition-number regimes. NAG-GS remains stable for relatively large values of $\alpha$ in these tests and often reaches the target tolerance in fewer iterations. The Heavy-Ball baseline uses momentum $0.9$, so its parameters are not matched to those of NAG-GS in \eqref{eq:hb-coeff-minimal}. The comparison therefore concerns these different parameter choices. Since \(\gamma=\mu\), for \(\kappa=10\) the theory gives
\[
    \alpha_*\simeq0.925,
    \qquad
    \alpha_{\rm stab}=1.
\]
For \(\kappa=1000\), the two values are both close to \(6.5\times10^{-2}\). Thus the optimal step approaches the stability boundary as the condition number grows, which is also visible in the narrow large-step region of the second panel.

\subsection{Nonsmooth strongly convex classification}
\label{sec:svm}

  We  train a binary linear support vector machine on the MNIST digits 4 and 9. The objective is
\begin{equation}
\label{eq:svm-objective}
    F(w)
    =
    \frac{\mu}{2}\|w\|^2
    +
    \frac{1}{m}\sum_{i=1}^m
    \bigl[1-y_i a_i^\top w\bigr]_+,
\end{equation}
where $[t]_+=\max\{0,t\}$. The hinge loss makes \eqref{eq:svm-objective} nonsmooth, while the quadratic term makes it $\mu$-strongly convex.

We use the standard MNIST train/test split restricted to digits 4 and 9, which gives 11,791 training samples and 1,991 test samples. Pixel values are scaled to $[0,1]$. We append a constant feature for the bias term and normalize every augmented feature vector so that $\|a_i\|\leq1$. We set $\mu=0.05$, use mini-batches of size 32, and train for 40 epochs. All results are averaged over five random seeds. For a given seed, all methods use the same shuffled sequence of mini-batches. We compute a high-accuracy deterministic SVM solution once and use its objective value, denoted by $F_{\rm ref}$, as a reference.

For NAG-GS, we use   the parameter schedule  of \Cref{thm:stoch-rate}: $\gamma_k\equiv\mu$ and
\[
    \alpha_k=\frac{2}{k+2}.
\]
No parameter of this schedule is tuned. We compare with three standard nonsmooth stochastic methods. The first is stochastic subgradient descent (SSGD) with step size $\eta_k=2/(\mu(k+2))$, projected onto the ball of radius $1/\mu$. This ball contains the minimizer: at optimality, $0=\mu w^*+h^*$ for a hinge-loss subgradient $h^*$ with $\|h^*\|\leq1$, hence $\|w^*\|\leq1/\mu$. The second is Pegasos \cite{shalevshwartz2011pegasos}, with $\eta_k=1/(\mu(k+1))$ and its standard projection step. The third is AdaGrad \cite{duchi2011adaptive}. For AdaGrad, the initial scale is selected from $\{0.1,0.3,1\}$ using a short five-epoch run with the first random seed; this gives $0.1$.

\begin{figure}[ht!]
    \centering
    \begin{subfigure}{0.48\linewidth}
        \includegraphics[width=\linewidth]{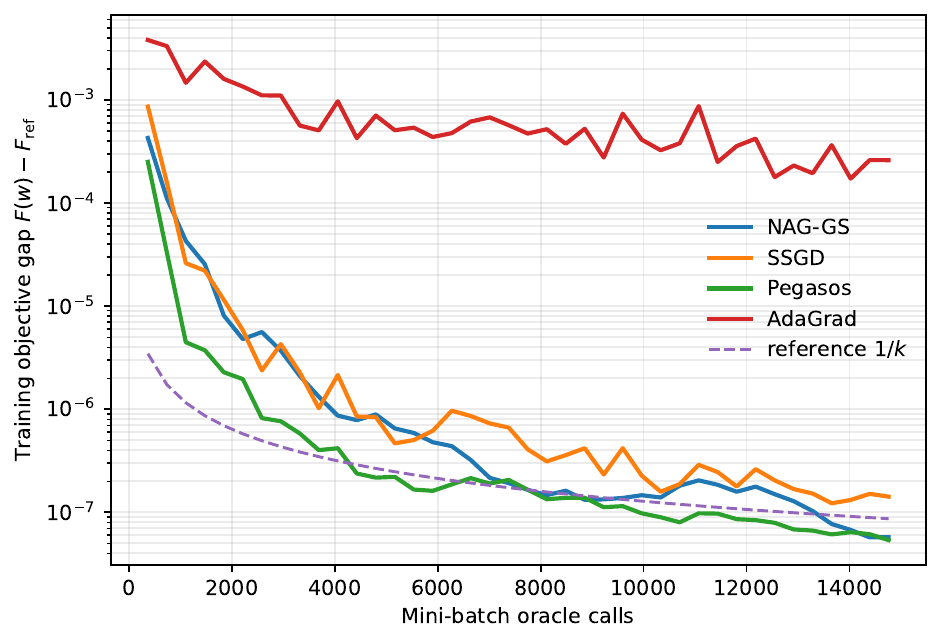}
        \caption{Training objective gap.}
    \end{subfigure}
    \hfill
    \begin{subfigure}{0.48\linewidth}
        \includegraphics[width=\linewidth]{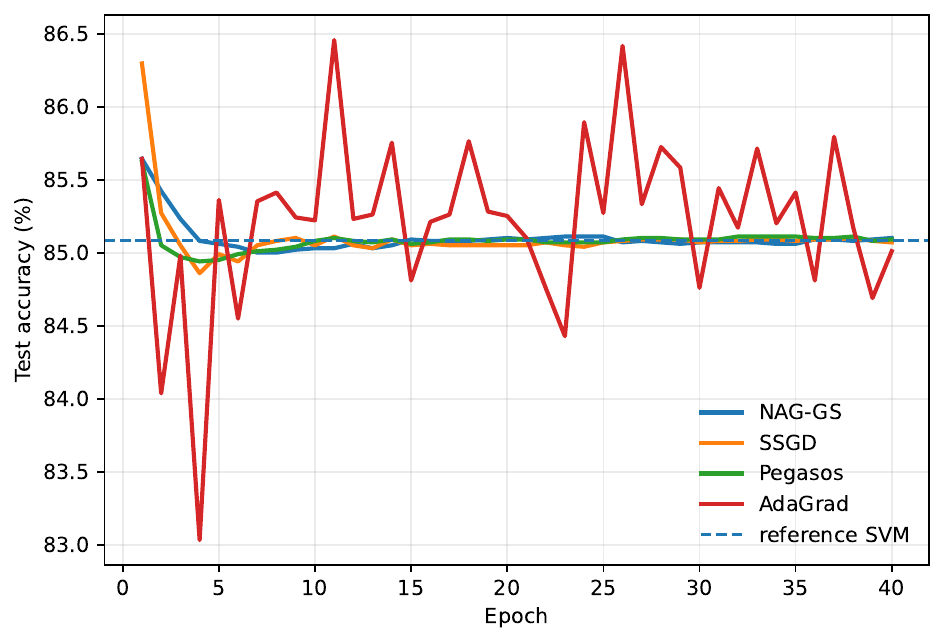}
        \caption{Test accuracy.}
    \end{subfigure}
    \caption{Nonsmooth strongly convex binary MNIST experiment. The curves show the mean over five random seeds. NAG-GS uses $\gamma=\mu$ and the decreasing step size from \Cref{thm:stoch-rate}.}
    \label{fig:svm-theory}
\end{figure}

 \Cref{fig:svm-theory} shows the optimization and classification results. NAG-GS rapidly approaches the reference objective and behaves similarly to SSGD and Pegasos. Pegasos gives the smallest mean objective gap at the end of the run, while NAG-GS gives a very similar value. SSGD also converges close to the reference solution, whereas AdaGrad is clearly slower for the chosen training budget. The test accuracies of NAG-GS, Pegasos and SSGD are essentially identical. NAG-GS has the largest mean test accuracy, although the differences between the three methods are very small.

This experiment evaluates the NAG-GS schedule from \Cref{thm:stoch-rate} on a nonsmooth strongly convex learning problem. With this schedule and no parameter tuning, NAG-GS performs competitively with SSGD and Pegasos. The dashed $1/k$ curve in \Cref{fig:svm-theory} is included only as a visual reference for the rate; it is not fitted to the data.

\begin{table}[ht!]
\centering
\caption{Binary MNIST classification with the strongly convex hinge-loss objective \eqref{eq:svm-objective}. The values are averaged over five random seeds. The objective gap is measured with respect to the high-accuracy reference solution.}
\label{tab:svm-theory}
\begingroup\small\setlength{\tabcolsep}{4pt}
 \begin{tabular}{lccc}
\toprule
Method & Gap after 10 epochs & Gap after 40 epochs & Test accuracy after 40 epochs \\
\midrule
NAG-GS
& $1.32\times10^{-6}\pm6.54\times10^{-7}$
& $5.71\times10^{-8}\pm2.99\times10^{-8}$
& $\mathbf{85.10\pm0.03}$ \\
SSGD
& $1.02\times10^{-6}\pm4.37\times10^{-7}$
& $1.42\times10^{-7}\pm6.96\times10^{-8}$
& $85.07\pm0.02$ \\
Pegasos
& $4.00\times10^{-7}\pm1.79\times10^{-7}$
& $\mathbf{5.37\times10^{-8}\pm8.70\times10^{-9}}$
& $85.09\pm0.02$ \\
AdaGrad
& $5.07\times10^{-4}\pm1.12\times10^{-4}$
& $2.61\times10^{-4}\pm9.67\times10^{-5}$
& $85.01\pm0.71$ \\
\bottomrule
\end{tabular}
 \endgroup
 \end{table}

 \paragraph{Heavy-ball equivalence and dynamic damping.}
\begingroup
We next consider three nonsmooth strongly convex problems. The first two are hinge-loss classification on normalized MNIST 4/9, with $(m,n)=(11791,785)$, and on \texttt{a9a}, with $(m,n)=(32561,124)$. The third is synthetic absolute-loss regression, with objective
\[
F(x)=\frac{\mu}{2}\|x\|^2+\frac{1}{m}\sum_{i=1}^{m}|a_i^\top x-b_i|.
\]
For this regression problem, $(m,n)=(2000,200)$, $\mu=0.01$, and the minimizer is specified when generating the data.
\endgroup
 All runs start at zero and use $\alpha_k=2/(k+2)$, batch size $32$, and uniform sampling with replacement. We report last queried iterates, without averaging; the five seeds are paired across configurations, and $K$ denotes the run length. $F_{\rm ref}$ is exact for regression; the SVM primal--dual residuals are below $10^{-10}$. With $\|a_i\|\le1$, write $g=\mu x+h$, $\|h\|\le1$: the updates keep $x,v$ in the ball of radius $1/\mu$, giving \eqref{eq:oracle} with $G=2$.
\begin{figure}[htbp]
 \centering\includegraphics[width=\linewidth]{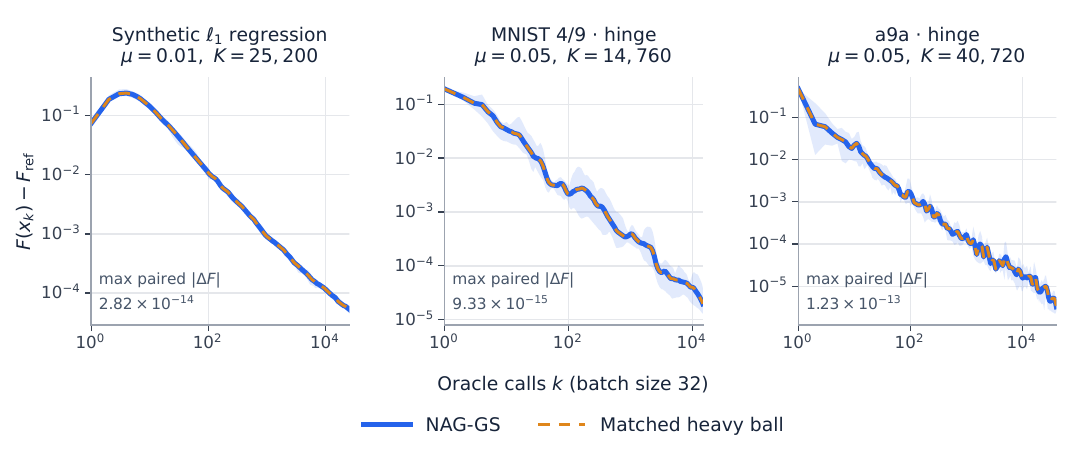}
 \caption{NAG-GS and matched HB ($\gamma=\mu$, $\alpha_k=2/(k+2)$), each evaluating its own subgradient on the same batches. Lines show means; shading shows the min--max range over five paired seeds. Annotations give the maximum paired objective difference across all seeds and recorded iterates.}
 \label{fig:hb-check}
\end{figure}
\begin{figure}[htbp]
 \centering\includegraphics[width=\linewidth]{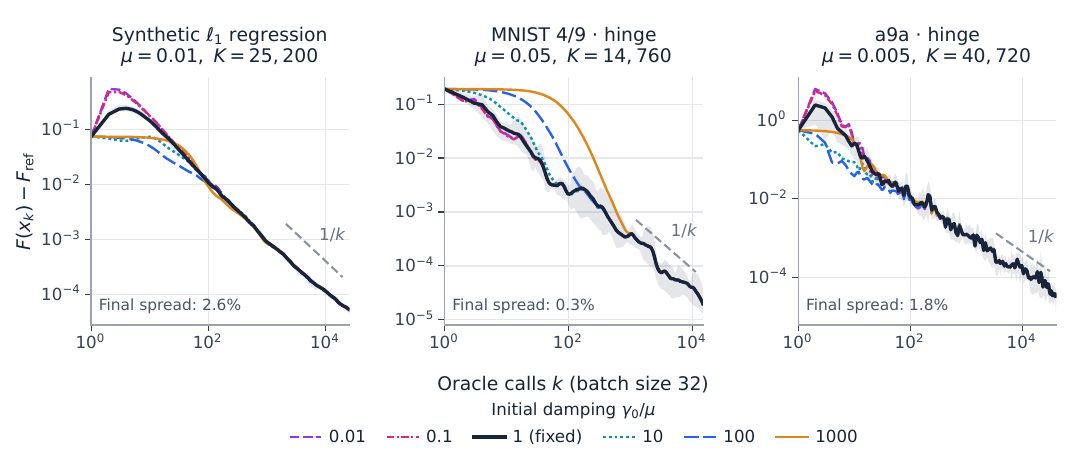}
 \caption{Dynamic update \eqref{eq:gamma-dynamic} with $\alpha_k=2/(k+2)$: means over five paired seeds. Shading shows the seed range for $\gamma_0=\mu$. Final spread compares the largest and smallest final mean gaps across the six initializations, relative to the smallest. Gray segments indicate a $1/k$ slope, not a bound.}
 \label{fig:gamma-check}
\end{figure}
\FloatBarrier
{\Cref{fig:hb-check} shows that NAG-GS and the matched Heavy-Ball update give nearly identical objective values, as predicted by their algebraic equivalence. \Cref{fig:gamma-check} shows the effect of choosing different initial damping values. The early iterations can differ, but the final objective gaps are close in these three tests.\par}

 \subsection{Logistic regression}
\label{sec:logistic}

We next train a multiclass logistic-regression model on MNIST \cite{lecun2010mnist}. We compare NAG-GS with SGD with momentum (SGD-M), AdamW, and AC-SA. The batch size is 500, the momentum in SGD-M is 0.9, and NAG-GS uses \(\mu=\gamma=1\). Weight decay is disabled for SGD-M and AdamW, so AdamW coincides with Adam in this experiment.

\Cref{tab:logreg} reports the best test accuracy for several learning rates. NAG-GS performs poorly at the smallest learning rate, but it remains effective at much larger values. With \(\alpha=0.5\), it reaches 92.38\% after 10 epochs. The other methods require smaller learning rates to obtain their best values in this comparison.

\begin{table}[ht!]
\centering
\caption{Multiclass logistic regression on MNIST. Accuracies are averaged over three random seeds.}
\label{tab:logreg}
\begin{tabular}{cccc}
\toprule
Learning rate & Epochs & Optimizer & Best test accuracy \\
\midrule
\multirow{4}{*}{$10^{-3}$} & \multirow{4}{*}{40}
& NAG-GS & $88.47\pm0.05$\\
&& SGD-M & $92.01\pm0.04$\\
&& AdamW & $\mathbf{92.71\pm0.01}$\\
&& AC-SA & $92.67\pm0.05$\\
\midrule
\multirow{4}{*}{$10^{-2}$} & \multirow{4}{*}{20}
& NAG-GS & $90.99\pm0.02$\\
&& SGD-M & $\mathbf{92.50\pm0.04}$\\
&& AdamW & $92.17\pm0.15$\\
&& AC-SA & $92.15\pm0.16$\\
\midrule
\multirow{4}{*}{$0.1$} & \multirow{4}{*}{15}
& NAG-GS & $\mathbf{92.18\pm0.07}$\\
&& SGD-M & $92.03\pm0.08$\\
&& AdamW & $89.83\pm0.24$\\
&& AC-SA & $90.93\pm0.15$\\
\midrule
\multirow{4}{*}{$0.5$} & \multirow{4}{*}{10}
& NAG-GS & $\mathbf{92.38\pm0.02}$\\
&& SGD-M & $90.04\pm0.35$\\
&& AdamW & $88.66\pm0.99$\\
&& AC-SA & $88.92\pm0.21$\\
\bottomrule
\end{tabular}
\end{table}

\begin{figure}[ht!]
\centering
\begin{subfigure}{0.43\textwidth}
\includegraphics[width=\linewidth]{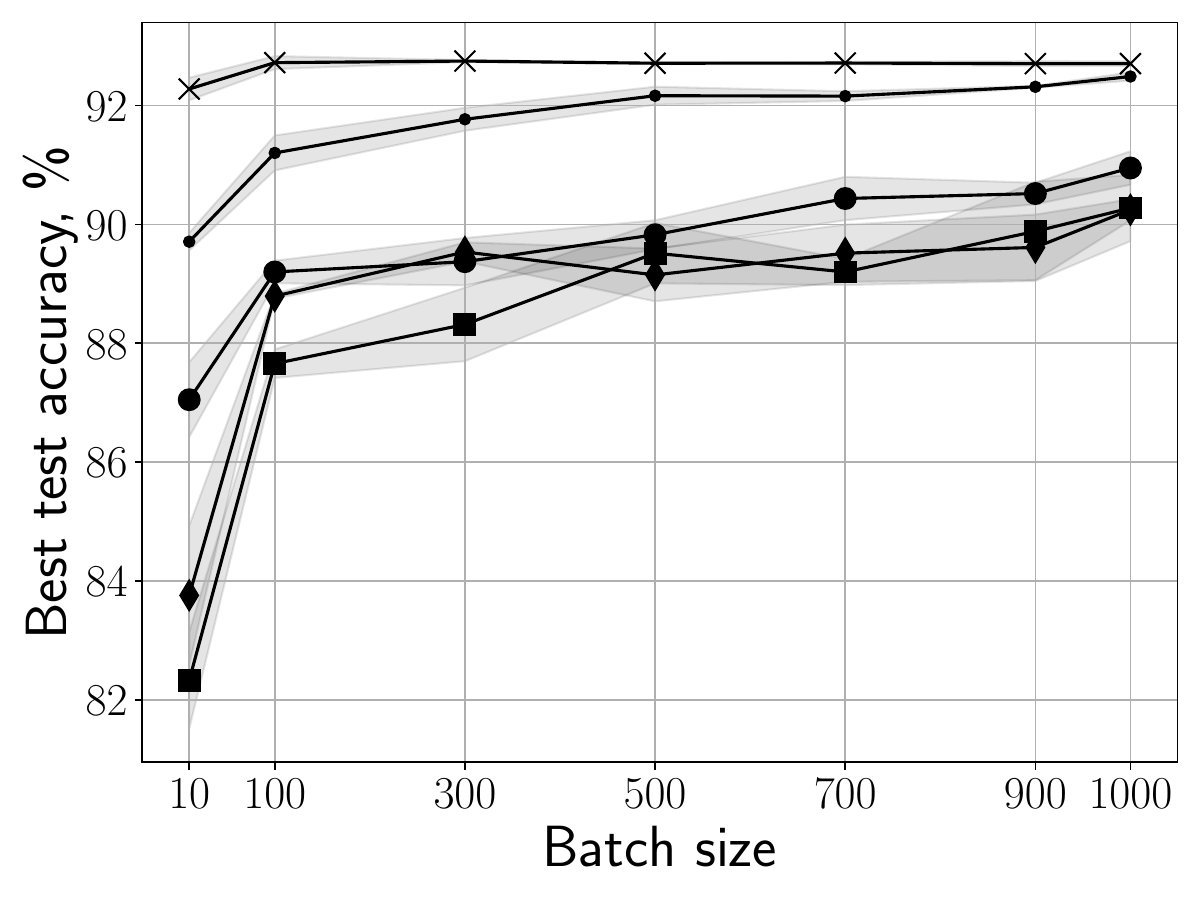}
\caption{AdamW}
\end{subfigure}
\hfill
\begin{subfigure}{0.43\textwidth}
\includegraphics[width=\linewidth]{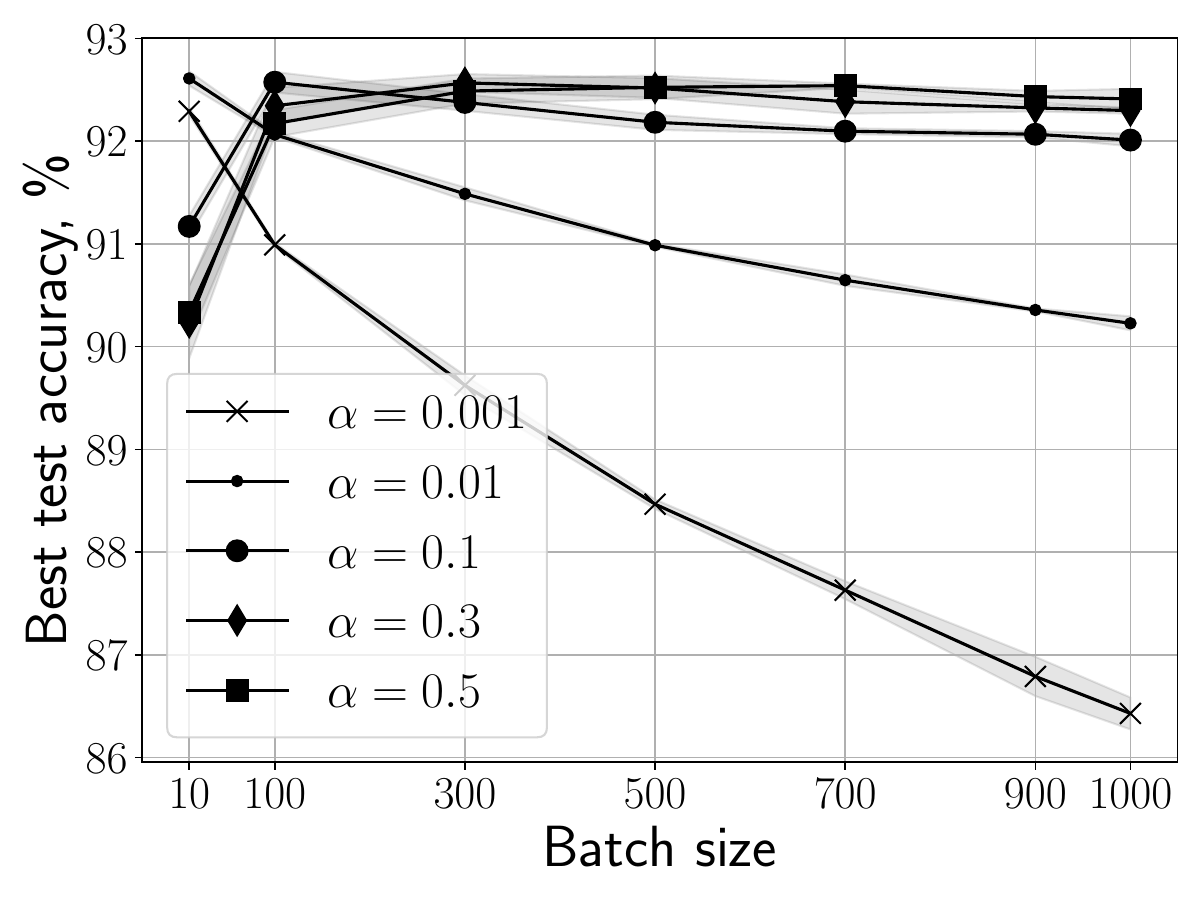}
\caption{NAG-GS}
\end{subfigure}
\\
\begin{subfigure}{0.43\textwidth}
\includegraphics[width=\linewidth]{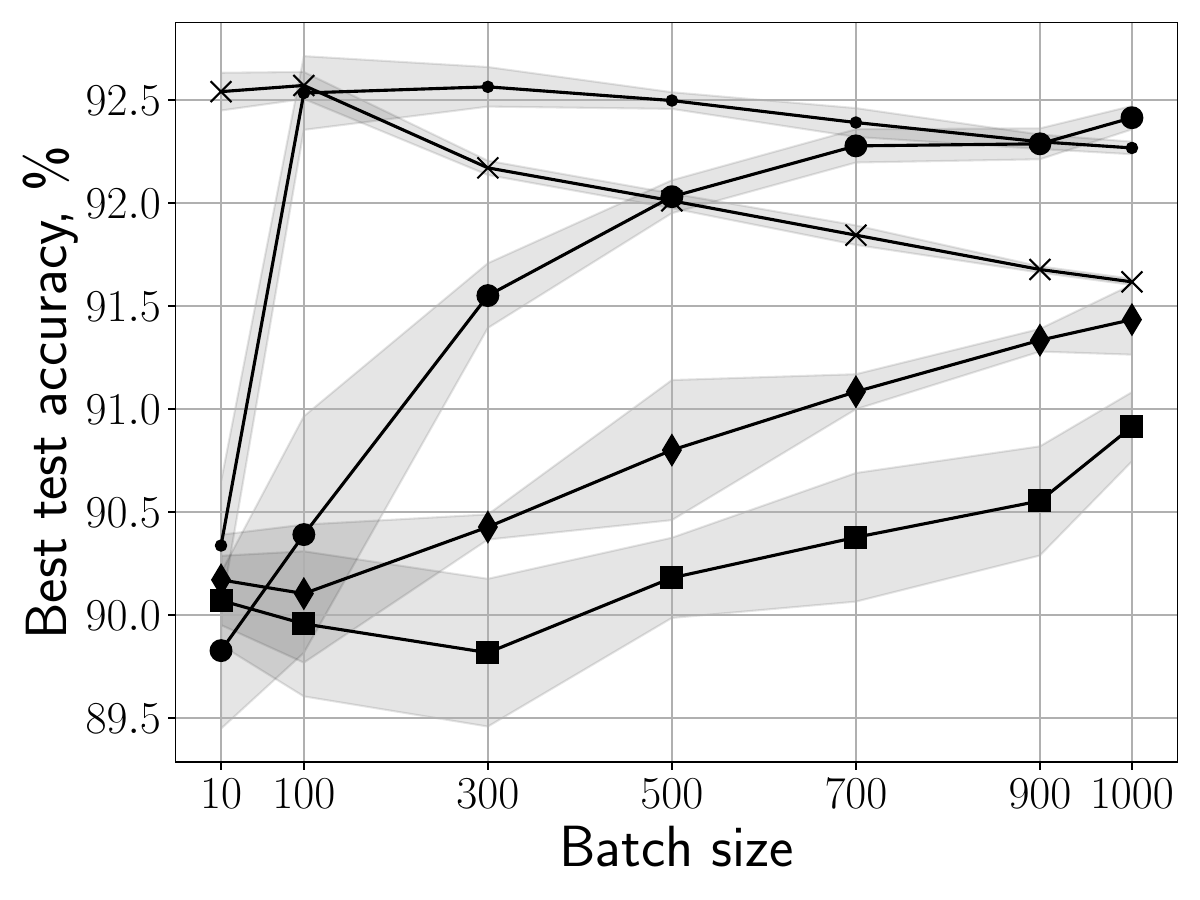}
\caption{SGD-M}
\end{subfigure}
\hfill
\begin{subfigure}{0.43\textwidth}
\includegraphics[width=\linewidth]{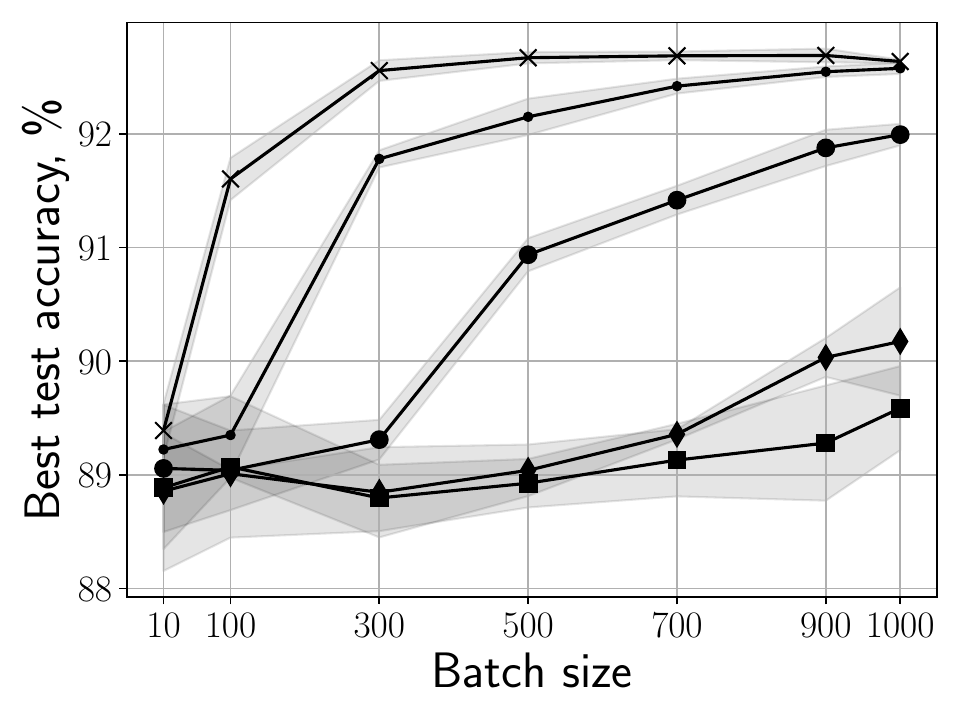}
\caption{AC-SA}
\end{subfigure}
\caption{Best MNIST test accuracy as a function of batch size for several learning rates. The learning-rate legend in the NAG-GS panel applies to all four plots.}
\label{fig:mnist-batch}
\end{figure}

We also vary the batch size and the learning rate. \Cref{fig:mnist-batch} shows the best test accuracy. The main observation is again the large-learning-rate regime: NAG-GS remains accurate for combinations for which the other optimizers degrade. At small learning rates, NAG-GS is not the best method.

This experiment is not meant to verify the \(O(1/k)\) theorem from \Cref{sec:stochastic}: the training protocol uses fixed learning rates and does not impose all assumptions of that theorem. Its role is instead to test whether the large-step behavior seen on quadratics persists on a simple stochastic convex model. The answer is positive in the tested range, but the advantage appears mainly at the larger learning rates.

\subsection{VGG-11 and ResNet-20}
\label{sec:cnn}

We now move to nonconvex stochastic problems. The goal changes at this point: the theory does not predict convergence rates for these models, so we only ask whether NAG-GS remains usable and competitive at comparatively large learning rates. We compare NAG-GS with SGD with momentum and weight decay (SGD-MW). We do not include AC-SA in these larger nonconvex experiments. Extending the comparison would require a separate tuning campaign, while AC-SA is already included in the convex tests above.

\paragraph{VGG-11.}
We train VGG-11 \cite{simonyan2014very} on CIFAR-10 \cite{krizhevsky2009learning} for 50 epochs with batch size 100. NAG-GS uses fixed \(\gamma=1\) and \(\mu=10^{-4}\). SGD-MW uses momentum 0.9 and weight decay \(10^{-4}\).

\begin{figure}[ht!]
\centering
\includegraphics[width=0.65\linewidth]{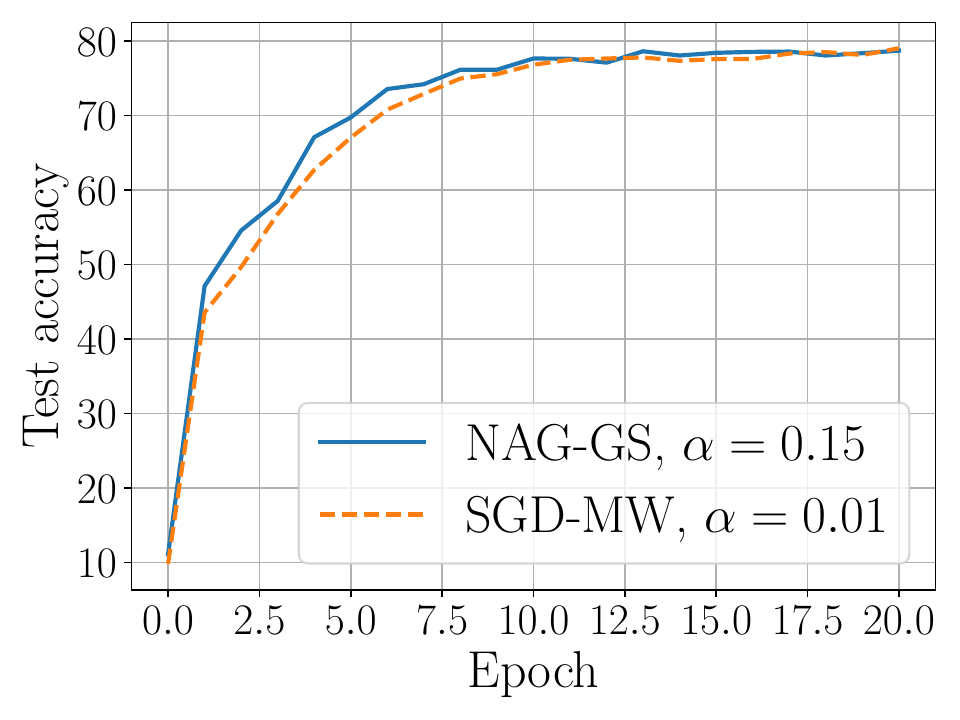}
\caption{VGG-11 on CIFAR-10. NAG-GS with $\alpha=0.15$ reaches high test accuracy faster in the first epochs than SGD-MW with $\alpha=0.01$. Curves are averaged over three random seeds.}
\label{fig:vgg-convergence}
\end{figure}

\Cref{tab:vgg,fig:vgg-convergence} show that NAG-GS can use learning rates 0.1 and 0.15 in this setup, while SGD-MW diverges at these values. With \(\alpha=0.15\), NAG-GS is also faster during the first epochs than SGD-MW with \(\alpha=0.01\). Later in training, the accuracies become similar.

\begin{table}[ht!]
\centering
\caption{Best VGG-11 test accuracy on CIFAR-10 over 50 epochs, averaged over three random seeds. A dash indicates divergence.}
\label{tab:vgg}
\begin{tabular}{ccc}
\toprule
Learning rate & Optimizer & Best test accuracy\\
\midrule
\multirow{2}{*}{$10^{-3}$} & NAG-GS & $66.15\pm0.81$\\
& SGD-MW & $78.84\pm0.21$\\
\midrule
\multirow{2}{*}{$10^{-2}$} & NAG-GS & $78.60\pm0.40$\\
& SGD-MW & $81.58\pm0.23$\\
\midrule
\multirow{2}{*}{$0.1$} & NAG-GS & $81.50\pm0.14$\\
& SGD-MW & --\\
\midrule
\multirow{2}{*}{$0.15$} & NAG-GS & $\mathbf{81.75\pm0.06}$\\
& SGD-MW & --\\
\bottomrule
\end{tabular}
\end{table}

\begin{figure}[ht!]
\centering
\begin{subfigure}{0.48\textwidth}
\includegraphics[width=\linewidth]{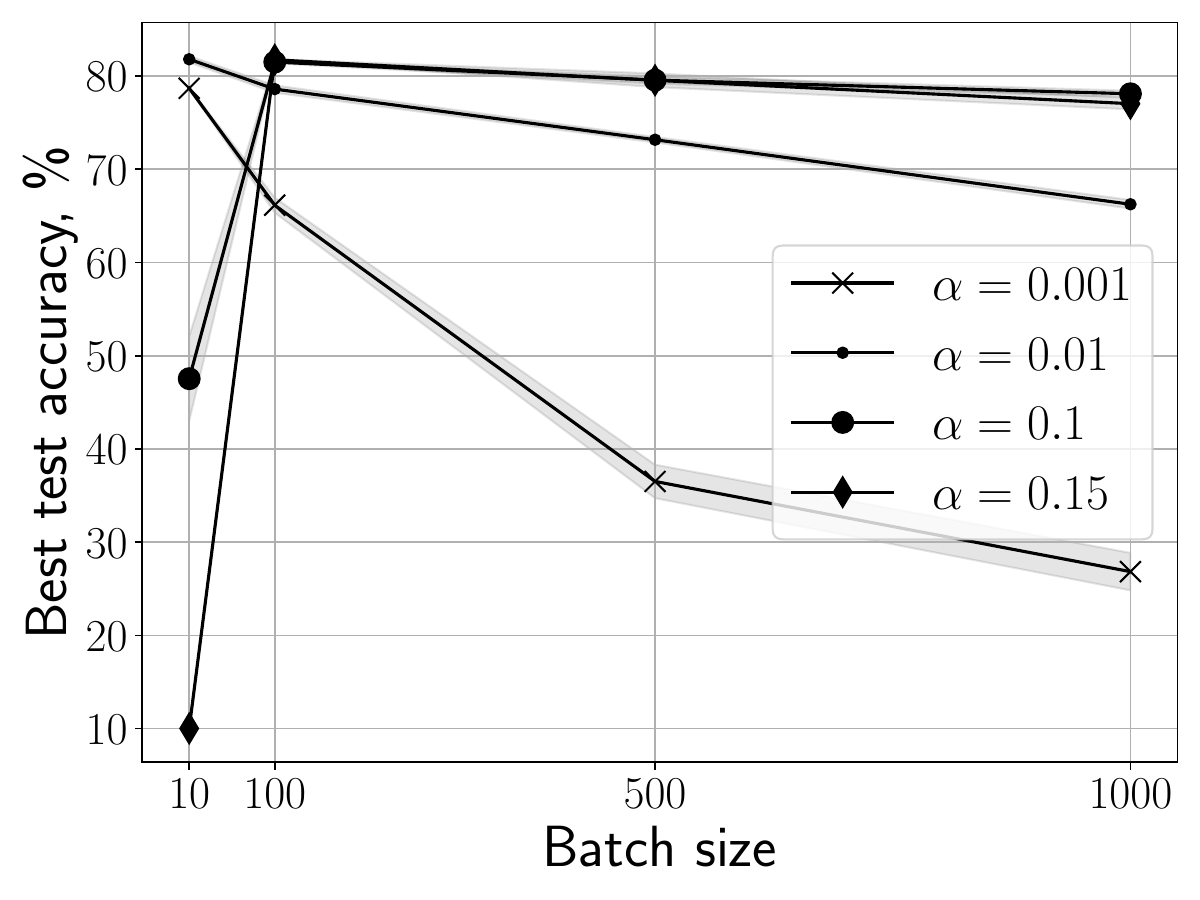}
\caption{NAG-GS}
\end{subfigure}
\hfill
\begin{subfigure}{0.48\textwidth}
\includegraphics[width=\linewidth]{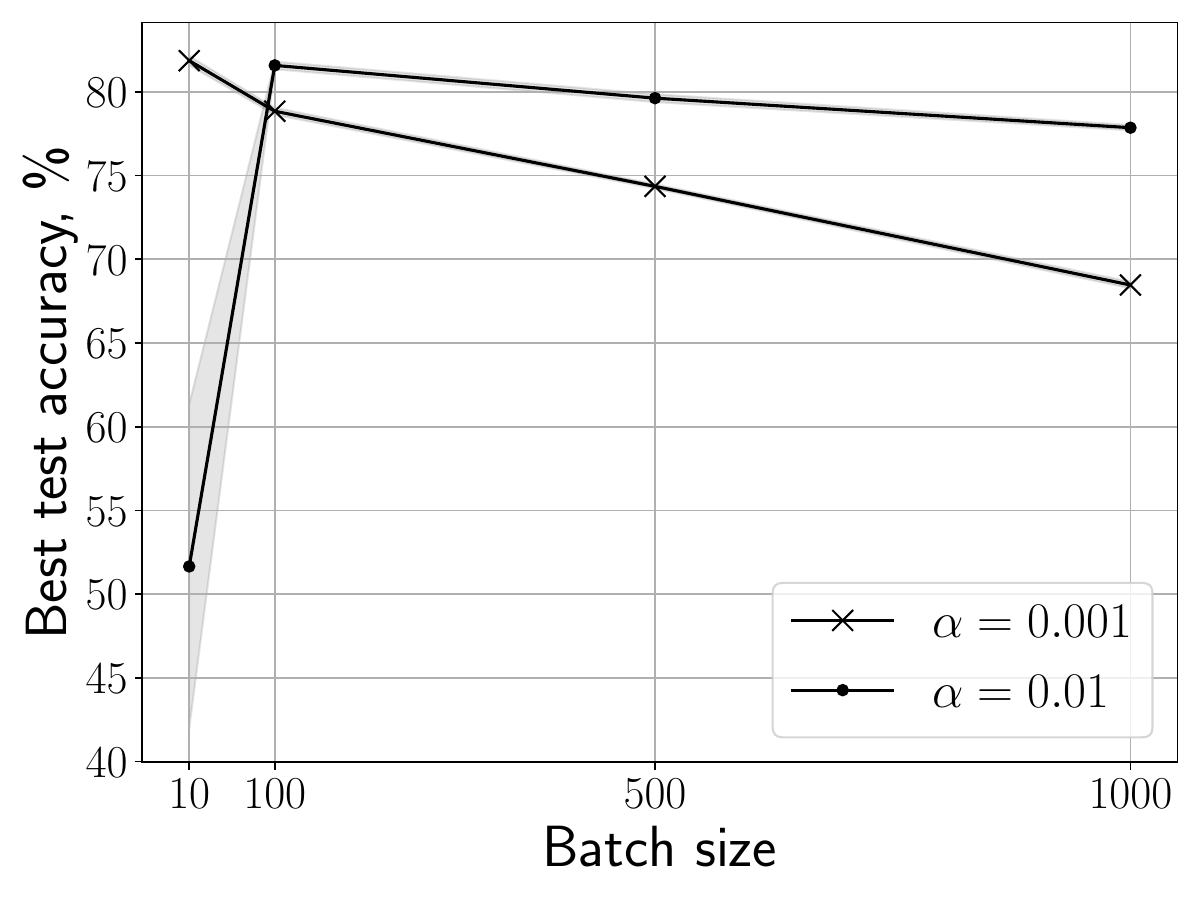}
\caption{SGD-MW}
\end{subfigure}
\caption{Best VGG-11 test accuracy as a function of batch size and learning rate.}
\label{fig:vgg-batch}
\end{figure}

The batch-size experiment in \Cref{fig:vgg-batch} gives the same qualitative picture. NAG-GS remains usable at the two largest learning rates. These learning rates lead to divergence or stagnation for SGD-MW and are therefore not shown in the SGD-MW panel. The main point is therefore robustness of the early training phase rather than a uniformly better final accuracy. At the smaller learning rates, SGD-MW can be as good as or better than NAG-GS.

\paragraph{ResNet-20.}
We train ResNet-20 \cite{he2016deep} on CIFAR-10. The best NAG-GS result in the search is obtained with \(\alpha=0.11\), \(\gamma=17\), and \(\mu=0.01\). The search space is
\[
\alpha\sim\operatorname{LogUniform}(10^{-2},10^2),\qquad
\gamma\sim\operatorname{LogUniform}(10^{-2},10^2),\qquad
\mu\sim\operatorname{Uniform}(-10,100),
\]
using Optuna \cite{akiba2019optuna}. Here \(\mu\) is a free algorithmic parameter. In particular, negative values in the search space are not related to the strong-convexity parameter used in the theory.

\begin{figure}[ht!]
\centering
\begin{subfigure}{0.49\textwidth}
\includegraphics[width=\linewidth]{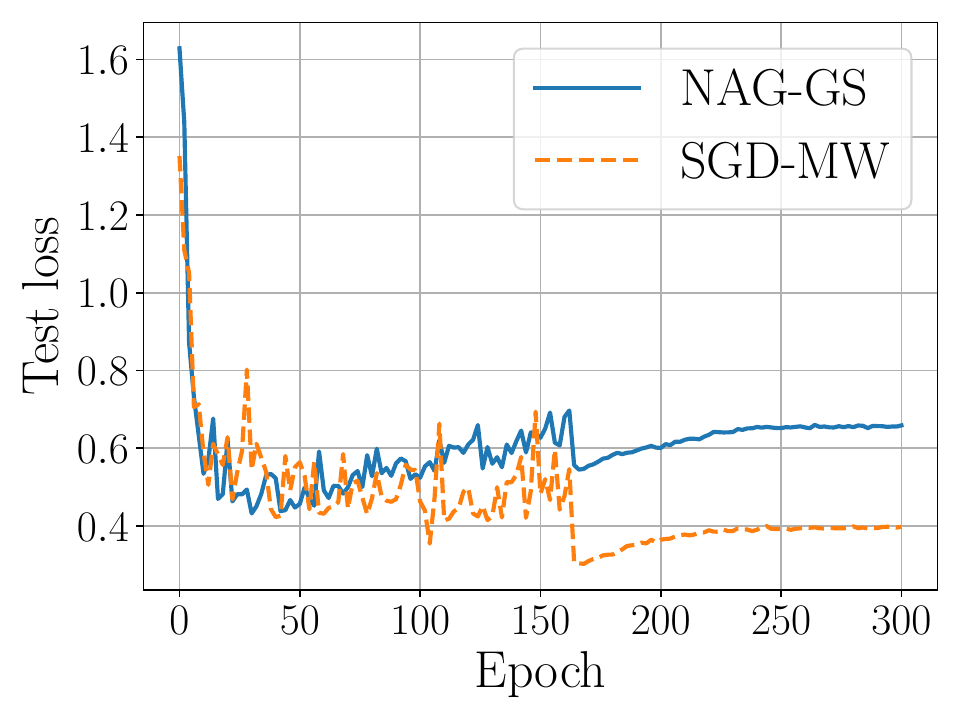}
\caption{Test loss}
\end{subfigure}
\hfill
\begin{subfigure}{0.49\textwidth}
\includegraphics[width=\linewidth]{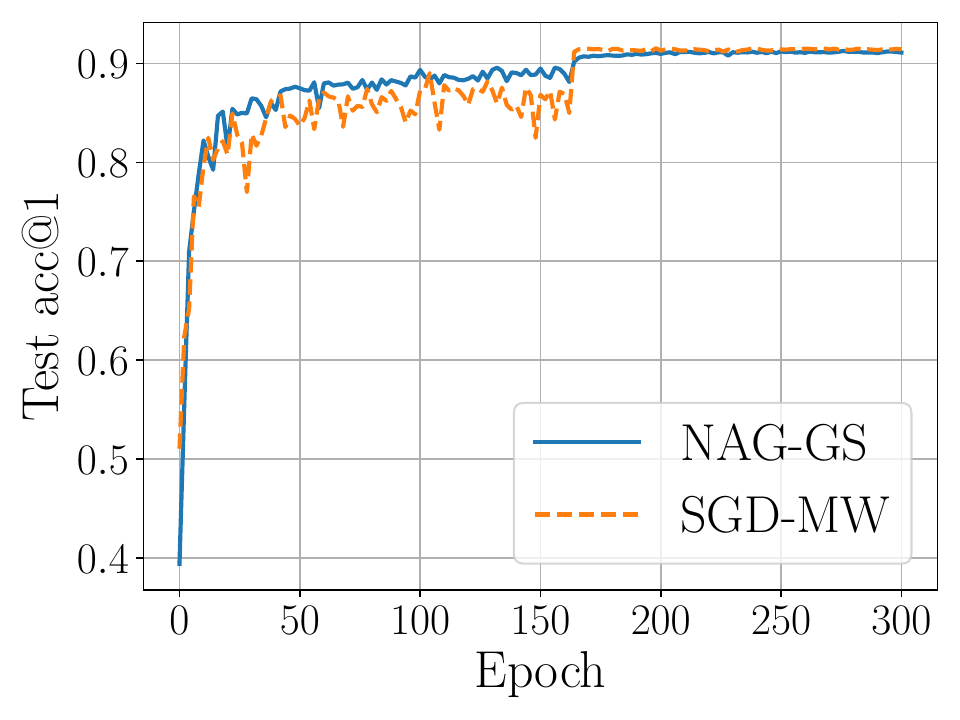}
\caption{Top-1 accuracy}
\end{subfigure}
\caption{NAG-GS and SGD-MW on ResNet-20/CIFAR-10.}
\label{fig:resnet}
\end{figure}

\Cref{fig:resnet} compares test loss and top-1 accuracy. NAG-GS has better test accuracy during the first 150 epochs, while the two methods become similar later in training. The hyperparameter search also shows that good configurations are not restricted to the line \(\gamma=\mu\) in this nonconvex problem; the corresponding phase diagrams are reported in \Cref{app:exp-extra}.

\subsection{Transformer models}
\label{sec:transformers}

\paragraph{RoBERTa.}
We fine-tune the pre-trained RoBERTa model \cite{liu2019roberta} using the \textsc{transformers} library \cite{wolf2020transformers} on the GLUE benchmark \cite{wang2018glue}. AdamW uses the reference polynomial learning-rate schedule. For NAG-GS, we tune a constant learning rate and fixed \(\gamma\), with \(\mu=1\). The learning-rate search is \([10^{-3},1]\), and the \(\gamma\) search is \([10^{-2},1]\).

\Cref{tab:glue} reports the best scores. NAG-GS is better on MRPC, tied on CoLA and WNLI, and slightly below AdamW on the remaining tasks. The selected learning rates are also on very different scales: around \(10^{-2}\) for NAG-GS and around \(10^{-5}\) for the reference AdamW setup. These values should not be compared as equivalent step sizes because the two optimizers use different update rules.

\begin{table}[ht!]
\centering
\caption{RoBERTa fine-tuning on GLUE.}
\label{tab:glue}
\resizebox{\textwidth}{!}{\begin{tabular}{lrrrrrrrrrr}
    \toprule
    Optimizer & \textsc{CoLA}  & \textsc{MNLI}  & \textsc{MRPC}  & \textsc{QNLI}  & \textsc{QQP}   & \textsc{RTE}   & \textsc{SST2}  & \textsc{STS-B} & \textsc{WNLI} \\
    \midrule
    AdamW     & \textbf{61.60} & \textbf{87.56} &         88.24  & \textbf{92.62} & \textbf{91.69} & \textbf{78.34} & \textbf{94.95} & \textbf{90.68} & \textbf{56.34} \\
    NAG-GS    & \textbf{61.60} &         87.24  & \textbf{90.69} &         92.59  &         91.01  &         77.97  &         94.50  &         90.21  & \textbf{56.34} \\
    \bottomrule
\end{tabular}
}
\end{table}

\paragraph{Vision Transformer.}
We fine-tune a Vision Transformer \cite{wu2020visual}, pre-trained on ImageNet \cite{deng2009imagenet}, on Food-101 \cite{bossard14}. We tune both NAG-GS and AdamW on a subset of the data using Weights \& Biases \cite{wandb}, then run the selected configurations on the full dataset. The NAG-GS search uses
\[
\alpha\sim\operatorname{LogUniform}(10^{-4},10^{-1}),\qquad
\mu\sim\operatorname{Uniform}(10^{-2},10),\qquad
\gamma\sim\operatorname{Uniform}(10^{-2},10).
\]
The selected NAG-GS parameters are \(\alpha=0.07929\), \(\gamma=0.3554\), \(\mu=0.1301\). The selected AdamW parameters are \(\mathrm{lr}=4.949\times10^{-5}\), \(\beta_1=0.8679\), and \(\beta_2=0.9969\).

\begin{table}[ht!]
\centering
\caption{Vision Transformer fine-tuning on Food-101.}
\label{tab:vit}
\begin{tabular}{ccc}
\toprule
Stage & NAG-GS & AdamW\\
\midrule
After 1 epoch & \textbf{0.8419} & 0.8269\\
After 25 epochs & \textbf{0.8606} & 0.8324\\
\bottomrule
\end{tabular}
\end{table}

In this experiment, NAG-GS gives the higher accuracy both after one epoch and after 25 epochs. Together with the RoBERTa results, this shows that the method can be used beyond convolutional architectures. We do not interpret the tuned values of \(\mu\) and \(\gamma\) through the strongly convex theory; they are simply algorithmic parameters in these tests.

\section{Discussion and conclusion}
\label{sec:conclusion}

NAG-GS is a simple semi-implicit momentum method with a large stable step-size range on quadratics. The analysis in \Cref{sec:quadratic} makes this statement precise. For fixed damping, the stability boundary follows from a two-dimensional Schur test, and the best fixed step can be written in closed form. With \(\gamma=\mu\), the asymptotic factor is the classical accelerated factor \((\sqrt\kappa-1)/(\sqrt\kappa+1)\).

The dynamic damping rule can also be kept in the theory. For a constant step, \(\gamma_k\) approaches \(\mu\) geometrically, and the resulting time-varying quadratic iteration inherits the stability of its limiting fixed-damping matrix. This does not give a complete nonlinear convergence theory for dynamic damping. It does give the dynamic damping rule a direct convergence guarantee on deterministic quadratic problems, even when \(\gamma_0\neq\mu\).

{For general strongly convex, possibly nonsmooth objectives, we use a time-varying Lyapunov function to establish an $O(1/k)$ convergence rate in expectation for the last iterate. The result holds with decreasing steps and dynamic damping from any $\gamma_0>0$, under a conditionally unbiased stochastic first-order oracle with bounded second moment.}

{The nonsmooth experiments illustrate the stochastic result and the effect of dynamic damping. They also confirm numerically the agreement with matched Heavy Ball.} The remaining numerical experiments concern nonconvex learning problems, which lie outside the assumptions of our convergence theory. They show that NAG-GS can remain effective at comparatively large learning rates in the tested settings.

Several questions remain open. {In particular, it would be interesting to extend the constant-step linear convergence result for dynamic damping beyond quadratic objectives.} Finally, the nonconvex experiments suggest that \(\mu\) and \(\gamma\) can be useful as {independent algorithmic parameters}, but there is currently no theory explaining how to choose them in this regime. These are natural directions for further work.

  \printbibliography[heading=bibintoc]

\appendix

\section{Proofs for the fixed-damping quadratic analysis}
\label{app:fixed-proof}

\subsection{Derivation of the modal matrix}

From the first equation in \eqref{eq:gs-discretization},
\[
    \xi_{k+1}=\frac{\xi_k+\alpha\nu_k}{1+\alpha}.
\]
The second equation gives
\[
    (\gamma+\alpha\mu)\nu_{k+1}
    =\gamma\nu_k+\alpha(\mu-\lambda)\xi_{k+1}.
\]
Substitution gives \eqref{eq:Tlambda}. A direct calculation then gives \eqref{eq:trace-det}--\eqref{eq:trace}.

\subsection{Proof of the fixed-damping stability theorem}

The characteristic polynomial of \(T_\lambda\) is
\[
    p_\lambda(z)=z^2-t_\lambda z+d,
\]
where \(t_\lambda=\tr T_\lambda\) and \(d=\det T_\lambda\). Since \(\alpha,\gamma,\mu>0\), we have \(0<d<1\).

For a real quadratic \(z^2-tz+d\), the roots are strictly inside the unit disk if and only if
\[
    1-t+d>0,\qquad 1+t+d>0,\qquad 1-d>0.
\]
Using \eqref{eq:trace-det}--\eqref{eq:trace}, these three quantities are
\begin{align*}
1-t_\lambda+d
&=\frac{\alpha^2\lambda}{(1+\alpha)(\gamma+\alpha\mu)}>0,\\
1+t_\lambda+d
&=\frac{4\gamma+2\alpha(\gamma+\mu)+\alpha^2(2\mu-\lambda)}{(1+\alpha)(\gamma+\alpha\mu)},\\
1-d
&>0.
\end{align*}
Only the second condition can fail. Its numerator decreases with \(\lambda\), so the worst case is \(\lambda=L\). If \(L\leq2\mu\), every term in the numerator is nonnegative and the constant and linear terms are strictly positive. Hence every \(\alpha>0\) is stable.

Assume now \(L>2\mu\). The stability condition becomes
\[
    (L-2\mu)\alpha^2-2(\gamma+\mu)\alpha-4\gamma<0.
\]
This quadratic has one positive root, which is \eqref{eq:alpha-stab}. The condition therefore holds exactly for \(0<\alpha<\alpha_{\rm stab}(\gamma)\). At equality, \(1+t_L+d=0\), so \(p_L(-1)=0\). Thus the boundary is reached through the eigenvalue \(-1\).
\qed

\subsection{Proof of the optimal fixed-step theorem}

Let
\[
    r=r_\gamma(\alpha)=\frac{\gamma}{\gamma+\alpha\mu}.
\]
For \(\lambda=\mu\), the matrix \(T_\mu\) is upper triangular. Its eigenvalues are
\[
    \frac{1}{1+\alpha}
    \quad\text{and}\quad
    r.
\]
Because \(\gamma\geq\mu\), we have \(r\geq1/(1+\alpha)\). Hence
\begin{equation}
\label{eq:lower-r}
    R(\alpha,\gamma)\geq r.
\end{equation}

We now ask when every modal eigenvalue lies in the disk \(|z|\leq r\). Set \(z=rw\). The characteristic polynomial becomes a quadratic in \(w\). Applying the non-strict Jury conditions gives
\begin{align*}
1-\frac{t_\lambda}{r}+\frac{d}{r^2}
&=\frac{\alpha^2(\lambda-\mu)}{\gamma(1+\alpha)}\geq0,\\
1+\frac{t_\lambda}{r}+\frac{d}{r^2}
&=\frac{4\gamma+2\alpha(\gamma+\mu)-\alpha^2(\lambda-\mu)}{\gamma(1+\alpha)},\\
1-\frac{d}{r^2}
&=\frac{\alpha(\gamma-\mu)}{\gamma(1+\alpha)}\geq0.
\end{align*}
Again, only the middle condition is restrictive, and it is worst at \(\lambda=L\). It holds exactly when
\[
    (L-\mu)\alpha^2-2(\gamma+\mu)\alpha-4\gamma\leq0,
\]
that is, when \(0<\alpha\leq\alpha_*(\gamma)\). Together with \eqref{eq:lower-r}, this proves
\[
    R(\alpha,\gamma)=r_\gamma(\alpha),
    \qquad 0<\alpha\leq\alpha_*(\gamma).
\]
This quantity is strictly decreasing in \(\alpha\).

It remains to show that the spectral radius cannot decrease again after \(\alpha_*\). Let \(r_*=r_\gamma(\alpha_*)\) and consider the characteristic polynomial of the \(\lambda=L\) mode. Using the equation that defines \(\alpha_*\), a direct substitution gives
\[
    p_L(-r_*)<0,
    \qquad \alpha>\alpha_*.
\]
For completeness, the numerator after factorization is
\[
-\gamma(\alpha-\alpha_*)(\alpha_*\mu+2\gamma)
\bigl(\alpha\alpha_*\gamma+2\alpha\alpha_*\mu+2\alpha\gamma
+\alpha_*^2\mu+2\alpha_*\gamma\bigr),
\]
while the denominator is positive. Since \(p_L(z)\to+\infty\) as \(z\to-\infty\), the polynomial has a real root smaller than \(-r_*\). Hence \(R(\alpha,\gamma)>r_*\) for every \(\alpha>\alpha_*\). This proves the result.
\qed

\section{Proof for dynamic damping}
\label{app:dynamic-proof}

For general positive steps, the dynamic rule satisfies
\[
\gamma_{k+1}-\mu
=\frac{\gamma_k-\mu}{1+\alpha_k},
\]
and iteration gives \eqref{eq:gamma-product}. Since \eqref{eq:gamma-dynamic} is a convex combination of \(\gamma_k\) and \(\mu\), the sequence remains positive and stays between \(\gamma_0\) and \(\mu\).

If \(\sum_k\alpha_k=+\infty\), then
\[
\sum_{k=0}^{\infty}\log(1+\alpha_k)=+\infty.
\]
Indeed, if \(\alpha_k\geq1\) for infinitely many \(k\), the claim is immediate. Otherwise \(\alpha_k<1\) eventually and \(\log(1+\alpha_k)\geq\alpha_k/2\). Hence the product in \eqref{eq:gamma-product} tends to zero and \(\gamma_k\to\mu\).

For \Cref{thm:dynamic}, the step is constant. We then have \eqref{eq:gamma-exact}, so the convergence to \(\mu\) is geometric. In particular, \(\gamma_k>0\) for every \(k\).

After diagonalizing \(A\), let
\[
    \T(\gamma)
    =\diag\bigl(T_{\lambda_1}(\alpha,\gamma),\ldots,T_{\lambda_n}(\alpha,\gamma)\bigr).
\]
The dynamic iteration is
\[
    z_{k+1}=\T(\gamma_k)z_k.
\]
The map \(\gamma\mapsto\T(\gamma)\) is continuous for \(\gamma>0\), and \eqref{eq:gamma-exact} implies
\[
    \T(\gamma_k)\longrightarrow \T(\mu).
\]
By assumption, \(\rho(\T(\mu))<1\).

Fix any \(q\) with \(\rho(\T(\mu))<q<1\). Choose \(q_0\) such that
\[
    \rho(\T(\mu))<q_0<q.
\]
We use the standard fact that, for every square matrix \(M\) and every \(\varepsilon>0\), one can choose a vector norm whose induced matrix norm satisfies
\[
    \|M\|_*<\rho(M)+\varepsilon.
\]
Apply this fact to \(M=\T(\mu)\) with \(\varepsilon=q_0-\rho(\T(\mu))\). We obtain a vector norm \(\|\cdot\|_*\) whose induced matrix norm satisfies
\[
    \|\T(\mu)\|_*<q_0.
\]
By continuity, there exists \(K\) such that
\[
    \|\T(\gamma_k)\|_*\leq q,
    \qquad k\geq K.
\]
Therefore
\[
    \|z_k\|_*
    \leq q^{k-K}\|z_K\|_*,
    \qquad k\geq K.
\]
All norms are equivalent in finite dimension, and the finite product that maps \(z_0\) to \(z_K\) can be absorbed into a constant. This gives the claimed linear convergence in the Euclidean norm.
\qed

\section{Proof of the stochastic Lyapunov bound}
\label{app:stochastic-proof}

We prove \Cref{lem:stoch-one-step} directly in the Euclidean setting. Write
\[
    \tau=\alpha_k,
    \qquad
    x=x_k,
    \qquad
    v=v_k,
    \qquad
    y=x_{k+1},
    \qquad
    w=v_{k+1},
\]
and let \(x^*\) be the minimizer.   Write $\gamma=\gamma_k$, $d=\gamma+\tau\mu=(1+\tau)\gamma_{k+1}$. The update \eqref{eq:naggs-state} gives
 \begin{equation}
\label{eq:xy-relation}
    y=\frac{x+\tau v}{1+\tau},
    \qquad
    w=  \frac{\gamma v+\tau\mu y-\tau g_{k+1}}{d} .
\end{equation}
Let
\[
    s=\E[g_{k+1}\mid\F_k]\in\partial f(y).
\]
Since \(f\) is \(\mu\)-strongly convex,
\[
    f(y)-f(x^*)
    \leq
    \langle s,y-x^*\rangle
    -\frac{\mu}{2}\|y-x^*\|^2.
\]
Convexity also gives
\[
    f(y)-f(x)
    \leq \langle s,y-x\rangle
    =\tau\langle s,v-y\rangle,
\]
where we used \(y-x=\tau(v-y)\). Combining the two inequalities,
\begin{equation}
\label{eq:function-part}
(1+\tau)\bigl(f(y)-f(x^*)\bigr)-\bigl(f(x)-f(x^*)\bigr)
\leq
\tau\langle s,v-x^*\rangle
-\frac{\mu\tau}{2}\|y-x^*\|^2.
\end{equation}

Set
\[
    u=v-x^*,
    \qquad z=y-x^*.
\]
Using the second relation in \eqref{eq:xy-relation} and taking conditional expectation,
\begin{align*}
&(1+\tau)  \frac{\gamma_{k+1}}{2} \E[\|w-x^*\|^2\mid\F_k]
-  \frac{\gamma}{2} \|u\|^2\\
&=  \frac{1}{2d}
 \E\left[\left\| \gamma  u+\tau \mu  z-  \tau  g_{k+1}\right\|^2\middle|\F_k\right]
-  \frac{\gamma}{2} \|u\|^2.
\end{align*}
Expanding the square and using \(\E[g_{k+1}\mid\F_k]=s\) gives
\begin{align*}
&(1+\tau)  \frac{\gamma_{k+1}}{2} \E[\|w-x^*\|^2\mid\F_k]
-  \frac{\gamma}{2} \|u\|^2\\
&=-  \frac{\gamma\mu\tau}{2d} \|u\|^2
+  \frac{\gamma\mu\tau}{d} \langle u,z\rangle
+  \frac{\mu^2\tau^2}{2d} \|z\|^2\\
&\quad-  \frac{\gamma\tau}{d} \langle u,s\rangle
-  \frac{\mu\tau^2}{d} \langle z,s\rangle
+  \frac{\tau^2}{2d} \E[\|g_{k+1}\|^2\mid\F_k].
\end{align*}
Adding this expression to \eqref{eq:function-part}, the inner-product and norm terms collect into
\begin{align}
\label{eq:pre-young}
\frac{1}{A_k}\left(\E[\mathcal E_{k+1}\mid\F_k]-\mathcal E_k\right)
\leq{}&
  \frac{\mu\tau^2}{d} \langle s,u-z\rangle
-  \frac{\mu\gamma\tau}{2d} \|u-z\|^2\\
&+  \frac{\tau^2}{2d}
 \E[\|g_{k+1}\|^2\mid\F_k].\nonumber
\end{align}
For any vectors \(s,r\),
\[
      \frac{\mu\tau^2}{d} \langle s,r\rangle-  \frac{\mu\gamma\tau}{2d} \|r\|^2
    \leq   \frac{\mu\tau^3}{2\gamma d} \|s\|^2.
\]
Also, Jensen's inequality gives
\[
    \|s\|^2
    \leq \E[\|g_{k+1}\|^2\mid\F_k].
\]
Applying these two bounds to \eqref{eq:pre-young} gives
\[
\frac{1}{A_k}\left(\E[\mathcal E_{k+1}\mid\F_k]-\mathcal E_k\right)
\leq
  \frac{\tau^2}{2\gamma}
 \E[\|g_{k+1}\|^2\mid\F_k],
\]
which is \eqref{eq:conditional-energy}.
\qed

\section{Additional quadratic illustrations}
\label{app:quadratic-extra}

The stationary covariance can be written explicitly on each quadratic mode, but the resulting rational formulas are long and are not needed for the stability proof. The Schur test in \Cref{thm:schur} gives the stability boundary directly. We therefore use the covariance only after stability has been established and include the following numerical illustrations to show two useful behaviors.

For \(\mu>L/2\), there is no finite stability boundary in \Cref{thm:schur}. Large steps can remain stable, although the fastest mean contraction occurs at a finite step. For \(\mu<L/2\), the stationary covariance becomes large as the step approaches the Schur boundary, and no stationary covariance exists once stability is lost.

The illustrations use three-dimensional quadratic problems of the form \(f(x)=\tfrac12(x-ce)^\top A(x-ce)\) with \(c=5\), \(\gamma=\mu=1\), and noise level \(\sigma=1\). The matrix \(A\) is obtained from a diagonal matrix with prescribed eigenvalues by a random orthogonal change of basis. We simulate \(200{,}000\) independent initial points drawn from a Gaussian distribution. In the first case, \(L=1.9\), so \(\mu>L/2\), and the tested steps are \(\alpha_*/2\), \(\alpha_*\), \(2\alpha_*\), and \(10\alpha_*\), with \(\alpha_*\simeq5.29\). In the second case, \(L=3\), so \(\mu<L/2\), and the tested steps are \(\alpha_*/2\), \(\alpha_*\), \((\alpha_*+\alpha_{\rm stab})/2\), and \(0.98\alpha_{\rm stab}\), with \(\alpha_*\simeq2.73\) and \(\alpha_{\rm stab}\simeq4.83\).

\begin{figure}[ht!]
\centering
\begin{subfigure}{0.48\textwidth}
\includegraphics[width=\linewidth]{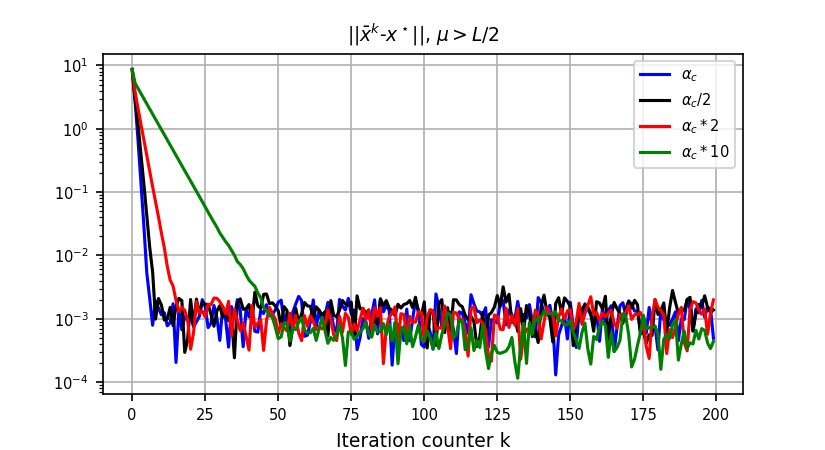}
\caption{$\mu>L/2$.}
\end{subfigure}
\hfill
\begin{subfigure}{0.48\textwidth}
\includegraphics[width=\linewidth]{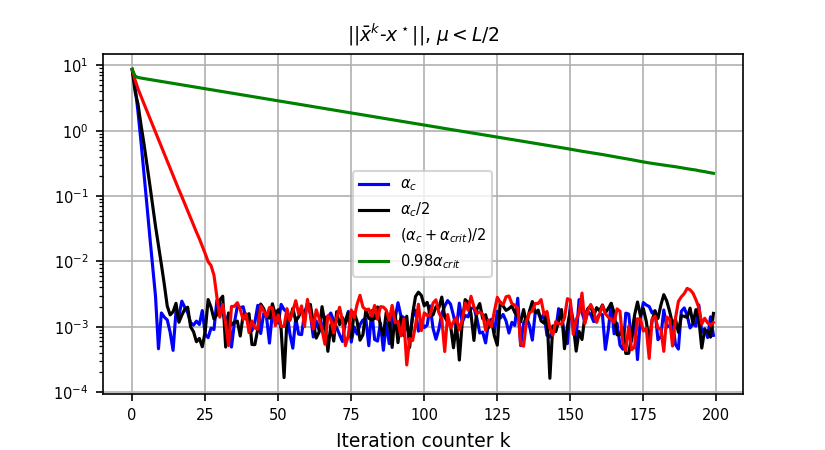}
\caption{$\mu<L/2$.}
\end{subfigure}
\caption{Distance between the empirical mean and the minimizer for several step sizes in the stochastic quadratic tests.}
\label{fig:quad-extra-rate}
\end{figure}

\begin{figure}[ht!]
\centering
\begin{subfigure}{0.32\textwidth}
\includegraphics[width=\linewidth]{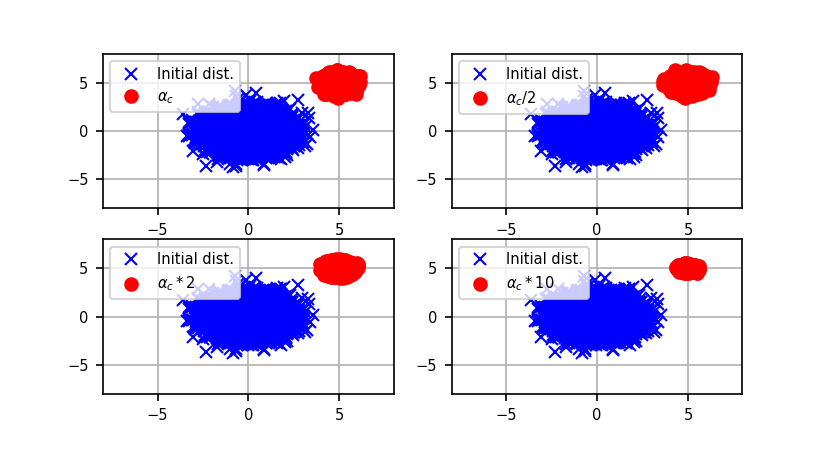}
\caption{$XY$.}
\end{subfigure}
\begin{subfigure}{0.32\textwidth}
\includegraphics[width=\linewidth]{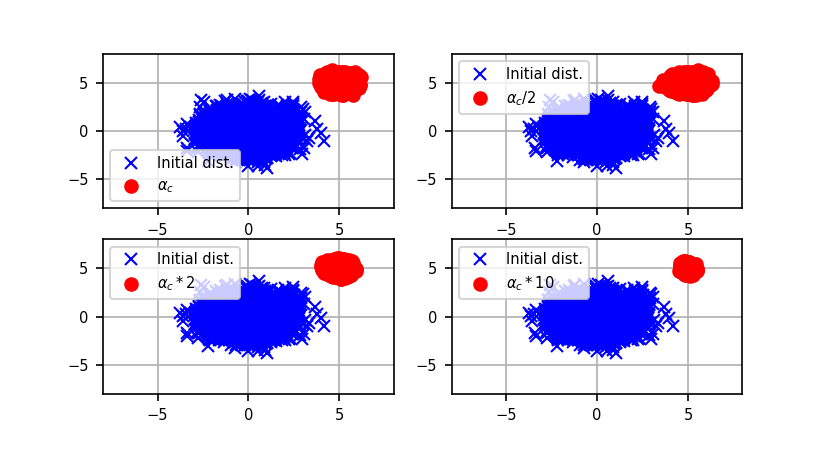}
\caption{$YZ$.}
\end{subfigure}
\begin{subfigure}{0.32\textwidth}
\includegraphics[width=\linewidth]{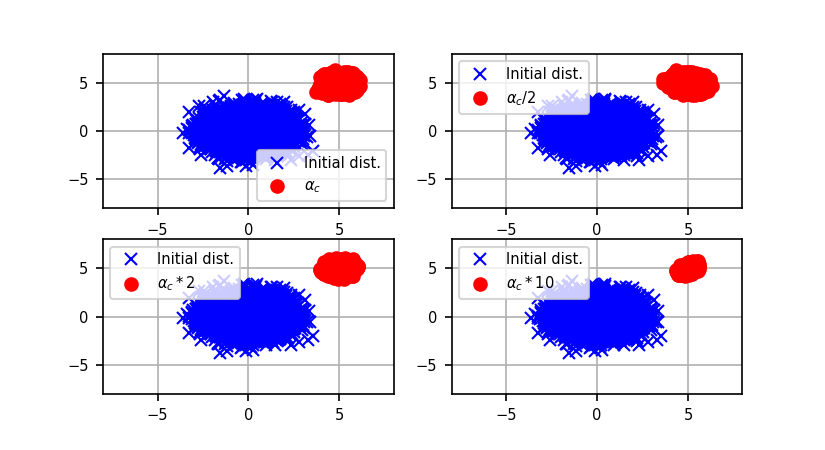}
\caption{$XZ$.}
\end{subfigure}
\caption{Initial and stationary samples in the case $\mu>L/2$.}
\label{fig:quad-extra-dist1}
\end{figure}

\begin{figure}[ht!]
\centering
\begin{subfigure}{0.32\textwidth}
\includegraphics[width=\linewidth]{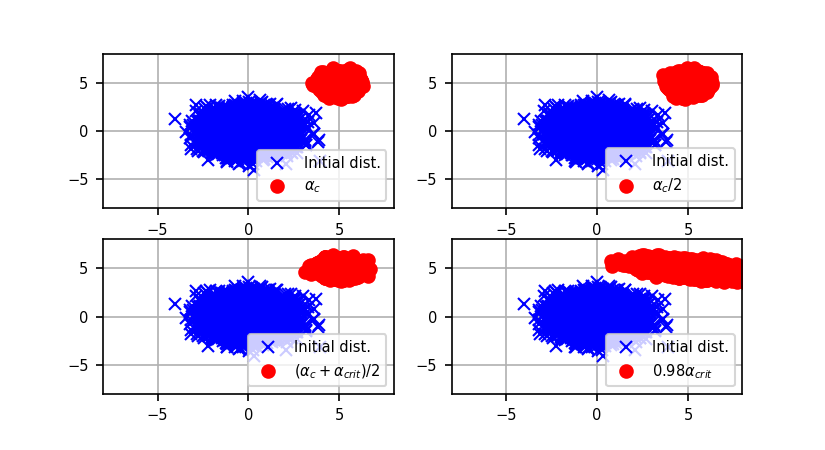}
\caption{$XY$.}
\end{subfigure}
\begin{subfigure}{0.32\textwidth}
\includegraphics[width=\linewidth]{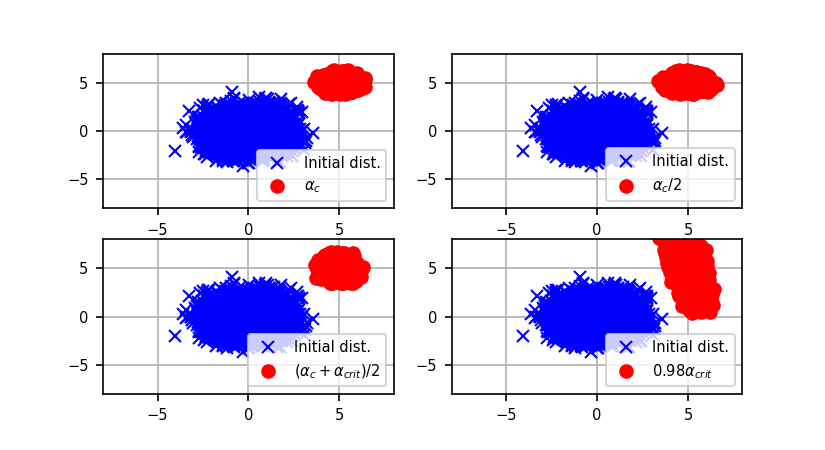}
\caption{$YZ$.}
\end{subfigure}
\begin{subfigure}{0.32\textwidth}
\includegraphics[width=\linewidth]{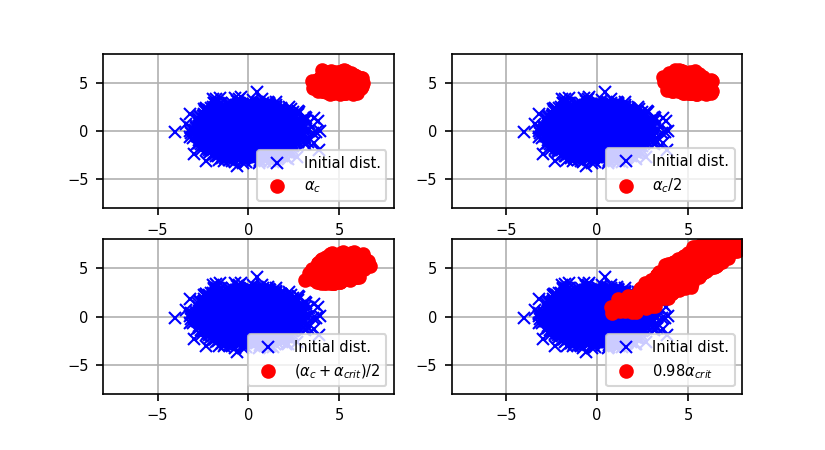}
\caption{$XZ$.}
\end{subfigure}
\caption{Initial and stationary samples in the case $\mu<L/2$. The spread increases when the step approaches the stability boundary.}
\label{fig:quad-extra-dist2}
\end{figure}

\section{Additional nonconvex SDE experiment}
\label{app:sde-extra}

We also include one small SDE experiment. Its purpose is only illustrative and it is not part of the convergence theory.

For the scalar stochastic gradient flow
\[
    \mathrm dx_t=-f'(x_t)\,\mathrm dt+\sigma\,\mathrm dW_t,
\]
the stationary density, when it exists and is normalizable, is proportional to
\[
    \exp\left(-\frac{2f(x)}{\sigma^2}\right).
\]
For the accelerated flow, we obtain empirical reference distributions by long Euler--Maruyama simulations. We then compare the empirical distributions using the Kullback--Leibler divergence, Wasserstein distance, and Kolmogorov--Smirnov statistic.

We use two scalar nonconvex functions,
\[
    f_1(x)=\frac{1}{50}\bigl(2\log(\cosh x)-5\bigr)^2,
    \qquad
    f_2(x)=\cos\left(1.6x+\frac{5}{3}\sin(0.64x)-\pi\right).
\]
For each function, we start from 100 points sampled uniformly over the plotted domain and integrate the stochastic dynamics independently. For \(f_1\), we use Euler--Maruyama step \(\alpha=8\times10^{-3}\), noise level \(\sigma=10^{-3}\), and \(\mu=1/33\). For \(f_2\), we use \(\alpha=1.5\), \(\sigma=10^{-2}\), and \(\mu=1\). These are finite-time comparisons: since \(f_2\) is bounded and periodic, the Gibbs density above is not normalizable on \(\mathbb R\), so this example does not establish convergence to a stationary probability law.

\begin{figure}[ht!]
\centering
\begin{subfigure}{0.90\textwidth}
\includegraphics[width=\linewidth]{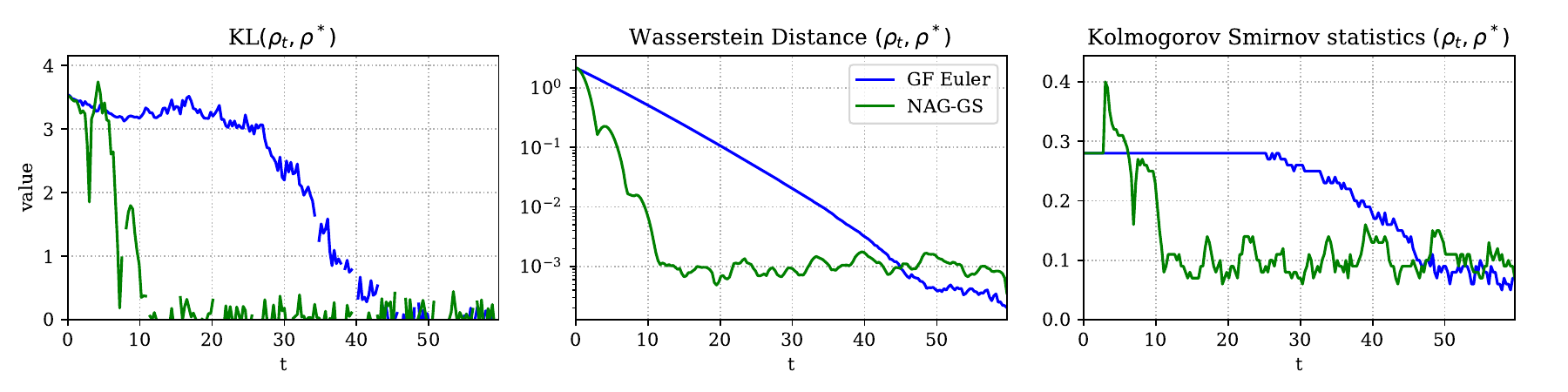}
\caption{Two-pits test function.}
\end{subfigure}
\hfill
\begin{subfigure}{0.90\textwidth}
\includegraphics[width=\linewidth]{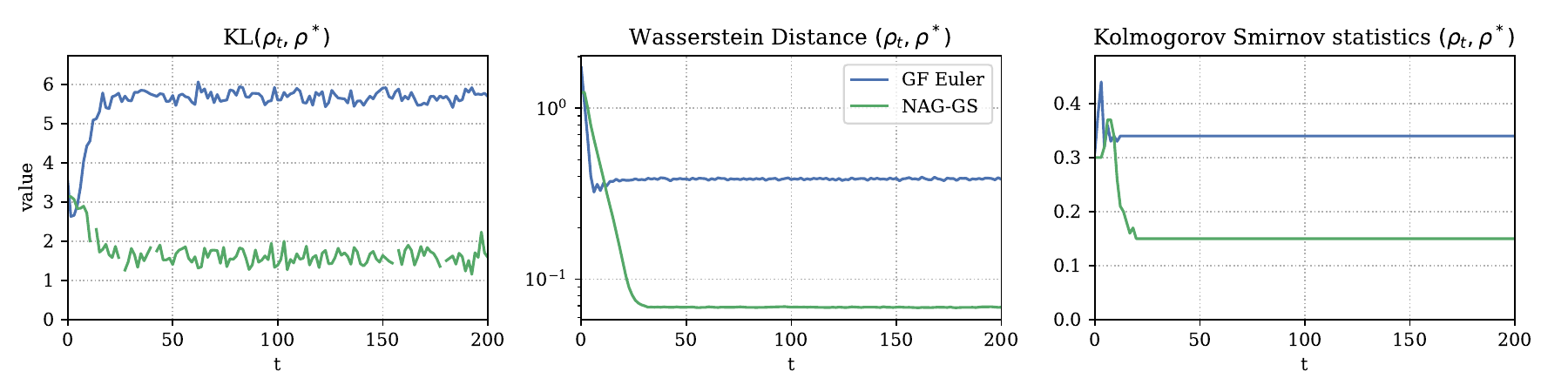}
\caption{Oscillatory test function.}
\end{subfigure}
\caption{Empirical distribution comparisons in two scalar nonconvex tests.}
\label{fig:sde-extra}
\end{figure}

\section{Additional experimental details}
\label{app:exp-extra}

Most experiments were run on a single NVIDIA V100 GPU from an NVIDIA DGX-1 system. The public implementation contains both PyTorch \cite{pytorch2017automatic} and JAX \cite{jax2018github,deepmind2020jax} versions of NAG-GS.

\subsection{Phase diagrams}

\begin{figure}[ht!]
\centering
\includegraphics[width=0.9\linewidth]{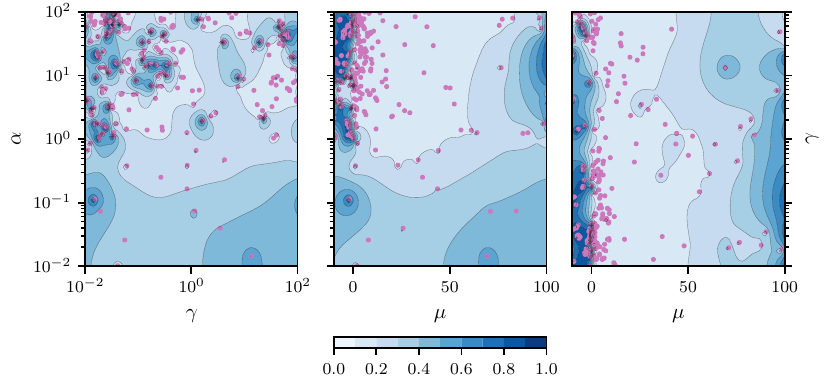}
\caption{Classification-error landscape for ResNet-20/CIFAR-10 under the NAG-GS hyperparameter search.}
\label{fig:phase-diagrams}
\end{figure}

The ResNet-20 hyperparameter search also gives a simple view of the empirical convergence regions. \Cref{fig:phase-diagrams} shows projections onto the \(\alpha\)-\(\gamma\), \(\alpha\)-\(\mu\), and \(\gamma\)-\(\mu\) planes. These plots are descriptive only; the nonconvex search is outside the scope of the quadratic theory.

\subsection{Optimizer-state and update cost}

NAG-GS stores one auxiliary vector. The measured update times from the implementation are reported in \Cref{tab:benchmark}. These timings measure optimizer overhead in the tested setup; they should not be read as hardware-independent complexity constants.

\begin{table}[ht!]
\centering
\caption{Measured duration of one optimizer update in the ResNet-20/CIFAR-10 implementation.}
\label{tab:benchmark}
\resizebox{\textwidth}{!}{\begin{tabular}{lrrrr}
    \toprule
    \textsc{Optimizer} & \textsc{Mean, s} & \textsc{Variance, s$^2$} & \textsc{Rel. Mean} & \textsc{Rel. Variance} \\
    \midrule
    SGD    & 0.458 & 0.008 &  1.0 &  1.0 \\
    NAG-GS & 1.648 & 0.045 &  3.6 &  5.5 \\
    SGD-M  & 3.374 & 0.042 &  7.4 &  5.2 \\
    SGD-MW & 3.512 & 0.037 &  7.7 &  4.7 \\
    AdamW  & 5.208 & 0.102 & 11.4 & 12.6 \\
    Adam   & 7.919 & 0.169 & 17.3 & 20.8 \\
    \bottomrule
\end{tabular}
}
\end{table}

\end{document}